\documentclass[journal]{IEEEtran}
\IEEEoverridecommandlockouts
\usepackage{amsmath,amsfonts,amssymb,euscript, graphicx, 
epsfig,
enumerate,float,afterpage, subfigure, ifthen}%

\newtheorem{thm}{Theorem}

\newtheorem{lem}{Lemma}

\newcommand{\expect}[1]{\mathbb{E}\left\{#1\right\}}
\newcommand{\defequiv}{\mbox{\raisebox{-.3ex}{$\overset{\vartriangle}{=}$}}}

\newcommand{\bv}[1]{{\boldsymbol{#1} }}
\newcommand{\script}[1]{{{\cal{#1} }}}

\begin{document}

\title
  {Universal Scheduling for Networks with Arbitrary Traffic, Channels, and Mobility}
\author{Michael J. Neely
\thanks{Michael J. Neely is with the  Electrical Engineering department at the University
of Southern California, Los Angeles, CA. (web: http://www-rcf.usc.edu/$\sim$mjneely).} 
\thanks{This material is supported in part  by one or more of 
the following: the DARPA IT-MANET program
grant W911NF-07-0028, 
the NSF Career grant CCF-0747525, and continuing through participation in the 
Network Science Collaborative Technology Alliance sponsored
by the U.S. Army Research Laboratory.}}

\markboth{}{Neely}

\maketitle

\begin{abstract} 
We extend stochastic network optimization theory to treat networks
with arbitrary sample paths for arrivals, channels, and mobility.  
The network can experience unexpected link or node failures, traffic bursts, and topology 
changes, and there are no probabilistic assumptions describing these time varying events. 
Performance of our scheduling algorithm is  
compared against an ideal $T$-slot lookahead policy that can 
make optimal decisions based on knowledge up to 
$T$-slots into the future. 
We develop a simple non-anticipating algorithm that provides network throughput-utility that is 
arbitrarily close to (or better than) that of the $T$-slot lookahead policy, with a tradeoff in the 
worst case queue backlog kept at any queue. 
The same policy offers even stronger performance, closely matching that
of an ideal \emph{infinite lookahead} policy, when ergodic assumptions are imposed. 
 Our analysis uses a sample path version of 
Lyapunov drift and 
provides a methodology for optimizing time averages in general time-varying
optimization problems. 
\end{abstract} 


\section{Introduction} 

\begin{figure}[htbp]
   \centering
   \includegraphics[height=1.5in, width=3in]{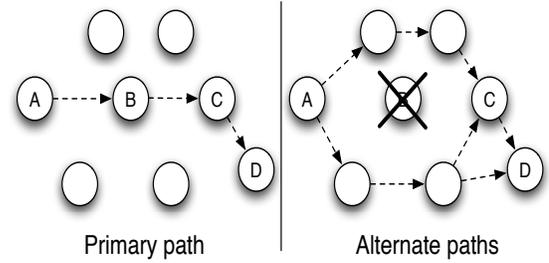} 
   \caption{A primary path from $A$ to $D$, 
   with alternative paths shown in the event of a failure at node $B$.}
   \label{fig:node-failure}
\end{figure}

Networks experience unexpected events.  Consider the network of Fig. \ref{fig:node-failure} and focus
on the session that sends a stream of packets from node $A$ to node $D$.  Suppose that several
paths are used, but due to congestion on other links, the primary path that can deliver the most data
is the path $A, B, C, D$.   However, suppose that there is a failure at node $B$ in the middle of the session.
An algorithm with perfect knowledge of the future would take advantage of the path $A, B, C, D$ while it is available, 
and would switch to alternate paths before the failure occurs.   The algorithm would also be able to predict 
the traffic load on different links at different times, and would optimally route in anticipation of these events. 

The above example holds if the network of Fig. \ref{fig:node-failure} is a wireline network, a wireless
network, or a mixture of wired and wireless connections.    
As another example, suppose the network contains an additional 
mobile wireless node $E$, and that  
the following unexpected event occurs:  Node $E$ moves into close
proximity to node $A$, allowing a large number of packets to be sent to it.  It then moves into close proximity to node $D$,
providing an opportunity to transmit packets to this destination node.  
If this event could be anticipated, we could take advantage of it and improve the short term throughput by routing many packets
over the relay $E$.  

These examples illustrate different types of unexpected events that can be exploited to improve performance.  There are of 
course even more complex sequences of arrival, channel, and mobility events that, if known in advance, could be exploited 
to yield improved performance. 
However, because realistic networks do not have knowledge of the future,  
it is not clear if these events can be practically used.  Surprisingly, this paper shows that it is possible to
reap the benefits of these time varying events without any knowledge of the future.  We show that a 
simple non-anticipating policy can closely track the 
performance
of an ideal \emph{$T$-slot lookahead policy} that has perfect knowledge of the future up to $T$ slots. 
Proximity to the performance of the $T$-slot lookahead policy comes with 
a corresponding tradeoff in the worst
case queue backlog stored in any queue of the network, which also affects a tradeoff in network delay. 

Specifically, we treat networks with slotted time with normalized slots $t \in \{0, 1, 2, \ldots\}$. We measure network  utility 
over an interval of timeslots according to a concave function of the time average throughput vector achieved 
over that interval.   We show that for any positive integer frame size $T$, 
and any interval that consists of $R$ frames of $T$ slots, 
the utility achieved over the interval is greater than or equal to the utility achieved by using the $T$-slot lookahead policy 
over each of the $R$ frames, minus a ``fudge factor'' that has the form: 
\[ \emph{fudge factor} = \frac{B_1T}{V}+ \frac{B_2V}{RT} \]
where $B_1$ and $B_2$ are constants, and $V$ is a positive parameter that can be chosen 
as desired to make the term 
$B_1T/V$ arbitrarily small, with a tradeoff in the worst case queue backlog that is $O(V)$.  
This shows that we reap almost the same
benefits of knowing the future up to $T$ slots if we choose $V$ suitably
large and if we wait for the completion of $R$ frames of size $T$, where $R$ is sufficiently large to make
$B_2V/(RT)$ small.   
Remarkably, the constants $B_1$ and $B_2$ can be explicitly 
computed in advance, without any assumptions on the underlying stochastic processes that describe the time varying events. 

This establishes a \emph{universal scheduling paradigm} that shows a 
single network algorithm can provide strong mathematical guarantees for
any network and for any time varying sample paths.  The algorithm that we use is not new:  It is a modified version 
of the backpressure based ``drift-minus-reward'' algorithms that we previously developed and used in different contexts 
in our prior work \cite{now}\cite{neely-fairness-infocom05}\cite{neely-fairness-ton}\cite{neely-thesis}. These algorithms were originally 
developed for 
the case when new arrivals and new channel states are independent and identically distributed (i.i.d.) 
over slots, and were analyzed using
a Lyapunov drift defined as an expectation over the underlying probability distribution.    However, it is known that
Lyapunov based algorithms that are designed under i.i.d. assumptions
yield similar performance under more general ergodic (but non-i.i.d.) probability models.
This is shown for network stability using a $T$-slot Lyapunov drift in \cite{tass-one-hop}\cite{neely-power-network-jsac}, and using a 
related delayed-queue analysis (that often provides tighter delay bounds) in \cite{neely-maximal-bursty-ton}.   Further, such algorithms
are known to be robust to non-ergodic situations such as when traffic yields ``instantaneous rates'' that can vary arbitrarily inside the capacity 
region \cite{neely-switch}\cite{neely-thesis}\cite{now}, and when ``instantaneous capacity regions'' can vary arbitrarily but are assumed
to always contain the traffic rate vector \cite{neely-mesh}.  However, the prior non-ergodic analysis \cite{neely-switch}\cite{neely-thesis}\cite{now}\cite{neely-mesh} still assumes an underlying probability 
model, and makes assumptions about traffic rates and network capacity with respect to this model. 

The analysis in this paper is new and uses a sample path version of Lyapunov drift, without any 
probabilistic assumptions. 
This framework allows treatment of realistic channels and traffic traces.
Because arbitrary  sample paths 
may not have well defined time averages, typical \emph{equilibrium} notions
of \emph{network capacity} and \emph{optimal time average utility} cannot be used.  We thus use a new metric
that  measures
performance with respect to ideal  $T$-slot lookahead policies.  
This is a possible framework for treating the 
important open questions identified in \cite{nequit-paper} 
concerning  \emph{non-equilibrium network theory}.  Further, our results provide universal techniques
for optimizing time averages that are useful for other types of time-varying systems. 

\subsection{Comparison to Related Work} 
We note that \emph{universal algorithms} are important in other fields.  For example, the 
universal Lempel-Ziv data compression algorithm operates on arbitrary files \cite{lempel-ziv},  
and universal stock portfolio 
allocation algorithms hold for arbitrary price sample paths \cite{cover-universal}\cite{universal-stock2}\cite{merhav-universal}\cite{neely-stock-arxiv}.
Our work provides a 
universal approach to network scheduling.  It is important to note that prior work in the area
of \emph{competitive ratio analysis} \cite{competitive-ratio-garay}\cite{plotkin-competitive}\cite{lin-shroff-large-N-energy}\cite{srikant-universal} and
\emph{adversarial queueing theory} \cite{andrews-max-profit}
also considers network scheduling problems with arbitrary sample paths, albeit in a different 
context. 
Work in \cite{plotkin-competitive} considers a large class of admission control problems for networks with 
random arrivals that earn revenue if accepted.  An algorithm is developed that yields revenue that
differs by a factor of $\Theta(\log(N))$ from that of an ideal algorithm with perfect knowledge of the future,
where $N$ is the number of network nodes.  Further, this asymptotic ratio is shown to be optimal, in the 
sense that there is always a \emph{worst case} sequence of packet arrivals that can reduce revenue by this
amount.   
Related $\Theta(\log(N))$ competitive ratio results are developed
for energy optimization in \cite{lin-shroff-large-N-energy} and for wireless
admission control in \cite{srikant-universal}.
The works \cite{competitive-ratio-garay}\cite{plotkin-competitive}\cite{lin-shroff-large-N-energy}\cite{srikant-universal} 
do not 
consider networks with time varying channels or mobility, and do not treat (or exploit)
network queueing.  An adversarial queueing theory 
example in \cite{andrews-max-profit} shows that, if channels are time varying, 
the competitive ratio can be much worse than logarithmic, even for a simple packet-based 
network with a single link.

Our work treats the difficult case of multi-hop networks with 
arbitrary time varying channels, mobility, and penalty constraints.  
However, rather than pursuing a competitive ratio analysis, we measure performance against
a $T$-slot lookahead metric.  We show that we can closely track the performance of an ideal $T$-slot
lookahead policy, for any arbitrary (but finite) $T$.  This does not imply that the algorithm has an optimal
competitive ratio, because the utility of a $T$-slot lookahead policy for finite $T$ may not be as good as the performance 
of an \emph{infinite} lookahead policy.  However, it turns out that our policy indeed approaches
an optimal competitive ratio (measured with respect to an infinite lookahead policy) under the special case
when the time varying events are ergodic.

Finally, we note that a frame-based metric, similar to our $T$-slot lookahead metric, 
was used in \cite{andrews-max-profit} to treat static wireline networks with arbitrary arrivals but
fixed topology and channel states.  There, an algorithm that queues all packets that arrive in a frame,
solves a network-wide utility maximization problem for these packets based on knowledge of the static 
link capacities, and implements the solution
in the next frame, is shown to achieve revenue that is close to that of a policy that a-priori knows the packet
arrival times over one frame.  In our context, we do not have the luxury of solving a network-wide
utility maximization problem based on known static link capacities, 
because the network itself is changing with time.    Our solution strategy 
is thus completely different from that of \cite{andrews-max-profit}.  Rather than a frame-based approach, 
our algorithm makes simple ``max-weight'' decisions every slot based on a quadratic Lyapunov function. 
It is interesting to note that this imposes a ``cost'' associated with each decision that depends on the current
network queue state, which is similar in spirit to the cost functions used in the algorithms
of \cite{competitive-ratio-garay}\cite{plotkin-competitive}\cite{lin-shroff-large-N-energy}\cite{srikant-universal} 
for competitive ratio analysis.

\subsection{Outline of Paper} 

The next section describes the general problem of constrained optimization 
of time varying systems.  
Section \ref{section:general-solution} provides a universal  solution
technique that measures performance against a $T$-slot lookahead policy. 
Section \ref{section:internet} applies the framework to a simple internet model,
and  Section \ref{section:network} applies the framework to a more extensive
class of time varying networks. 

\section{General Time Varying Optimization Problems} \label{section:general-framework}

Here we provide a
framework for universal 
constrained optimization for a general class of 
time varying systems.  The framework is applied in Sections \ref{section:internet} and
\ref{section:network} to solve
the network problems of interest. 

Consider a slotted system with normalized timeslots $t \in \{0, 1, 2, \ldots\}$. The system contains
$K$ queues with current backlog given by the vector $\bv{Q}(t) = (Q_1(t), \ldots, Q_K(t))$. 
Let $\omega(t)$ 
denote a \emph{random event} that occurs on slot $t$. The random event $\omega(t)$ represents
a collection of current system parameters and takes values in some abstract \emph{event 
space} $\Omega$.   
We treat $\omega(t)$ as a pre-determined function of $t \in \{0, 1, 2, \ldots\}$, 
although each value $\omega(t)$ is not revealed until the beginning of slot $t$. 
 Every slot, a 
system controller observes the current value of $\omega(t)$ and chooses a \emph{control action} $\alpha(t)$, 
constrained to some \emph{action space} $\script{A}_{\omega(t)} $ 
that can depend on $\omega(t)$.   The random event $\omega(t)$ and the corresponding control action $\alpha(t) \in \script{A}_{\omega(t)}$ produce a \emph{service vector} $\bv{b}(t) = (b_1(t), \ldots, b_K(t))$, an
\emph{arrival vector} $\bv{a}(t) = (a_1(t), \ldots, a_K(t))$, and two  
\emph{attribute vectors} $\bv{x}(t) = (x_1(t), \ldots, x_M(t))$, 
$\bv{y}(t) = (y_0(t), y_1(t), \ldots, y_L(t))$ (for some non-negative integers $M$ and $L$).
These vectors are general functions of $\omega(t)$ and $\alpha(t)$: 
\begin{eqnarray*}
a_k(t) &=& \hat{a}_k(\alpha(t), \omega(t)) \: \: \: \forall k \in \{1, \ldots, K\} \\
b_k(t) &=& \hat{b}_k(\alpha(t), \omega(t)) \: \: \: \forall k \in \{1, \ldots, K\} \\
x_m(t) &=& \hat{x}_m(\alpha(t), \omega(t)) \: \: \: \forall m \in \{1, \ldots, M\}  \\
y_l(t) &=& \hat{y}_l(\alpha(t), \omega(t)) \: \: \: \forall l \in \{0, 1, \ldots, L\} 
\end{eqnarray*}

The queue dynamics are determined by the arrival and service variables by: 
\begin{equation} \label{eq:q-update} 
Q_k(t+1) =  \max[Q_k(t) - b_k(t), 0] + a_k(t) 
\end{equation} 

Let $\overline{x}_m$ be the time average of $x_m(t)$ under a particular control policy
implemented over a finite number of slots $t_{end}$: 
\begin{equation} \label{eq:time-av-x}
 \overline{x}_m \defequiv \frac{1}{t_{end}}\sum_{\tau=0}^{t_{end}-1} x_m(\tau) 
 \end{equation} 
Define $\overline{a}_k$, $\overline{b}_k$, $\overline{y}_l$  similarly.  Define $\overline{\bv{x}} \defequiv (\overline{x}_1, \ldots, \overline{x}_M)$. 
The goal is to design a policy 
that solves the following time average optimization problem: 
\begin{eqnarray}
\mbox{Minimize:} &  \overline{y}_0 + f(\overline{\bv{x}}) \label{eq:opt} \\
\mbox{Subject to:} & \overline{y}_l + g_l(\overline{\bv{x}}) \leq 0 \: \: \: \forall l \in \{1, \ldots, L\} \label{eq:c1} \\
&\overline{a}_k \leq \overline{b}_k \: \: \: \forall k \in \{1, \ldots, K\} \label{eq:c2}\\
&\overline{\bv{x}} \in \script{X} \label{eq:c3} \\
& \alpha(t) \in \script{A}_{\omega(t)} \: \: \forall t \in \{0, \ldots, t_{end}-1\} \label{eq:c4}
\end{eqnarray}
where $f(\bv{x})$, $g_l(\bv{x})$ are convex  \emph{cost functions} 
of the vector  $\bv{x} = (x_1, \ldots, x_M)$, and $\script{X}$ is a general 
convex subset of $\mathbb{R}^M$. 
The above problem is of interest even if there are no underlying queues $Q_k(t)$ (so that the constraints
(\ref{eq:c2}) are removed), and/or if $\script{X} = \mathbb{R}^M$ (so that the constraint (\ref{eq:c3}) is removed), 
and/or if $f(\cdot) = g_l(\cdot) = 0$.

The above problem is stated in terms of a finite horizon of size $t_{end}$. 
The minimum cost in (\ref{eq:opt}) is defined for a given $\omega(t)$ function over $t \in \{0, \ldots, t_{end}-1\}$, 
and considers all possible control actions that can be implemented over the time horizon, 
including actions that have full knowledge of the future
values of $\omega(t)$.   However, we desire a \emph{non-anticipating} control policy that only knows
the current $\omega(t)$ value on each slot $t$, and has no knowledge of the future.   We note that 
a theory for solving stochastic network optimization 
problems similar to (\ref{eq:opt})-(\ref{eq:c4}) (without the set constraint (\ref{eq:c3}))
is developed in \cite{now} for an infinite horizon context under the assumption that $\omega(t)$ is i.i.d. over slots with 
some (possibly unknown) probability distribution, 
and related problems are treated in a fluid limit sense in \cite{stolyar-greedy}\cite{stolyar-gpd-gen}.
Here, 
we do not consider any probability model for $\omega(t)$, so that it is not possible to use law of large number
averaging principles, or to achieve the optimum
performance in an ``expected'' sense.  Rather, 
we shall deterministically achieve an ``approximate''
optimum for large $t_{end}$ values, to be made precise in future sections.

\subsection{Boundedness and Feasibility Assumptions} 

Here we present assumptions concerning the 
functions $\hat{a}_k(\cdot)$, $\hat{b}_k(\cdot)$, $\hat{x}_m(\cdot)$, $\hat{y}_l(\cdot)$
that ensure the problem (\ref{eq:opt})-(\ref{eq:c4}) is feasible with a bounded infimum cost. 

\emph{Assumption A1:} The functions $\hat{a}_k(\cdot)$, $\hat{b}_k(\cdot)$, 
$\hat{x}_m(\cdot)$, $\hat{y}_l(\cdot)$  
are bounded, so that  
for all $\omega \in \{\omega(0), \ldots, \omega(t_{end}-1)\}$, all  $\alpha \in  \script{A}_{\omega}$, 
and all $k \in \{1, \ldots, K\}$, $m \in \{1, \ldots, M\}$,  $l \in \{0, 1, \ldots, L\}$ we have: 
\begin{eqnarray*}
0 \leq \hat{a}_k(\alpha, \omega) \leq a_k^{max} \\
0 \leq \hat{b}_k(\alpha, \omega) \leq b_k^{max} \\
x_m^{min} \leq \hat{x}_m(\alpha, \omega) \leq x_m^{max} \\
y_l^{min} \leq \hat{y}_l(\alpha, \omega) \leq y_l^{max} 
\end{eqnarray*}
for some finite constants $a_k^{max}$, $b_k^{max}$, $x_m^{min}$, $x_m^{max}$,  $y_l^{min}$, 
$y_l^{max}$. 
Further, the cost functions $f(\bv{x})$ and $g_l(\bv{x})$ are defined over all vectors $(x_1, \ldots, x_M)$ 
that satisfy $x_m^{min} \leq x_m \leq x_m^{max}$ for all $m \in \{1, \ldots, M\}$, 
and have finite upper and lower bounds
$f^{min}$, $f^{max}$, $g_{l}^{min}$, $g_{l}^{max}$ over this region.

\emph{Assumption A2:} For all $\omega \in \{\omega(0), \ldots, \omega(t_{end}-1)\}$, there is at least one 
control action $\alpha_{\omega}' \in \script{A}_{\omega}$ that
satisfies: 
\begin{eqnarray} 
\hat{y}_l(\alpha_{\omega}', \omega) + g_l(\hat{\bv{x}}(\alpha_{\omega}', \omega)) \leq 0 &   \forall l \in \{1, \ldots, L\} \label{eq:a2-1} \\
\hat{a}_k(\alpha_{\omega}', \omega) \leq \hat{b}_k(\alpha_{\omega}', \omega) &  \forall k \in \{1, \ldots, K\}  \label{eq:a2-2} \\
\hat{\bv{x}}(\alpha_{\omega}', \omega) \in \script{X} \label{eq:a2-3} 
\end{eqnarray} 
where $\hat{\bv{x}}(\alpha, \omega)$ is defined: 
\[ \hat{\bv{x}}(\alpha, \omega) \defequiv (\hat{x}_1(\alpha, \omega), \ldots, \hat{x}_M(\alpha, \omega)) \]

For a given $\{\omega(0), \ldots, \omega(t_{end}-1)\}$, we say that the problem (\ref{eq:opt})-(\ref{eq:c4}) 
is \emph{feasible} 
if there exist control actions $\{\alpha(0), \ldots, \alpha(t_{end}-1)\}$ that satisfy the constraints (\ref{eq:c1})-(\ref{eq:c4}). 
Assumption A2 ensures that the problem is feasible (just
consider the control actions $\alpha(t) = \alpha_{\omega(t)}'$ for all $t$, and use Jensen's inequality to 
note that $g_l(\overline{\bv{x}}) \leq \overline{g_l(\bv{x})}$). 
Define $F^*$ as the infimum value of the cost metric (\ref{eq:opt}) over all feasible policies. The value $F^*$ is finite
by Assumption A1. 
Then for any $\epsilon>0$, there are control actions $\{\alpha^*(0), \ldots, \alpha^*(t_{end}-1)\}$
that satisfy the constraints (\ref{eq:c1})-(\ref{eq:c4}) with a total cost that satisfies: 
\[ F^* \leq \overline{y}_0^* + f(\overline{\bv{x}}^*) \leq F^* + \epsilon \]
Appendix B provides conditions that ensure the infimum $F^*$ is achievable by a particular policy (i.e., with $\epsilon = 0$). 
We note that Assumption A2 can be relaxed to only require feasibility over frames of $T$ slots, although 
 the resulting performance bounds in Theorem \ref{thm:1}a and \ref{thm:1}b are slightly altered in this case. 

\subsection{Cost Function Assumptions} 

The cost functions $f(\bv{x})$, $g_l(\bv{x})$ are assumed to be convex and continuous
over the region of all $(x_1, \ldots, x_M)$ vectors that satisfy: 
\begin{equation} \label{eq:x-interval} 
x_m^{min} \leq x_m \leq x_m^{max} \: \: \mbox{ for all 
$m \in \{1, \ldots, M\}$} 
\end{equation} 
In addition, we  
assume that the magnitude of the $m$th left and right partial derivatives of $f(\bv{x})$ with respect to $x_m$
are upper bounded by finite constants $\nu_m \geq 0$ for all $\bv{x}$ that satisfy (\ref{eq:x-interval}) 
and such that $x_m^{min} < x_m < x_m^{max}$.\footnote{All convex functions have well defined
right and left partial derivatives.}  Similarly, the magnitude of the right and left partial derivatives of $g_l(\bv{x})$ are
upper bounded by finite constants $\beta_{l,m} \geq 0$.  This implies that: 
\begin{eqnarray}
&f(\bv{x} + \bv{y}) \leq f(\bv{x}) + \sum_{m=1}^M\nu_m |y_m|  \label{eq:cost-assumption-f} \\
&g_l(\bv{x} + \bv{y}) \leq g_l(\bv{x})  + \sum_{m=1}^M\beta_{l,m}|y_m|  \label{eq:cost-assumption-g} 
\end{eqnarray}
for all $\bv{x}$, $\bv{y}$ such that $\bv{x}$ and $\bv{x} + \bv{y}$
are in the region specified by (\ref{eq:x-interval}).

An example of a non-differentiable cost function that satisfies all of the above assumptions is: 
\[ f(x_1, \ldots, x_M) = \max[x_1, \ldots, x_M] \]
In this case, we have $\nu_m = 1$ for all $m$.  Another example is a separable cost
function: 
\begin{equation} \label{eq:separable-structure} 
f(x_1, \ldots, x_M) = \sum_{m=1}^M f_m(x_m) 
\end{equation} 
where  functions $f_m(x)$ are continuous and convex with 
derivatives bounded in magnitude by $\nu_m$ over $x_m^{min}\leq x \leq x_m^{max}$.

\subsection{Applications of This Optimization Framework} 

We show in Section \ref{section:network} that this problem applies to general dynamic 
networks.  There, the  
$\omega(t)$ value represents a collection of channel conditions for all network links 
on slot $t$.  This includes the simple model where link conditions are either in the $ON$ or $OFF$ state, 
representing connections or disconnections that can vary from slot to slot due to fading channels and/or
user mobility.  Each node can discover the state of its links by probing to find existing neighbors on the
current slot.  The $\alpha(t)$
value represents a collection of routing, resource allocation, and/or flow control decisions that are
taken by the network in reaction to the current $\omega(t)$ value.  

In addition to dynamic networks and queueing systems, 
the problem (\ref{eq:opt}) has applications in many other areas that involve
optimization over time varying systems.  For example, our recent work in 
\cite{neely-stock-arxiv} presents 
applications to  stock market trading problems.
There, $\omega(t)$ represents a vector of current stock prices, and 
the constraints $\overline{a}_k \leq \overline{b}_k$ ensure that the average amount
of stock $k$ sales cannot exceed the average amount of stock $k$ purchases. 

We note that Assumption A2, which assumes that for any $\omega$ there exists 
an action $\alpha_{\omega}'$ that
satisfies $\hat{a}_k(\alpha_{\omega}', \omega) \leq \hat{b}_k(\alpha_{\omega}', \omega)$, 
often holds for systems that have a physical ``idle'' control action that reduces the inequality 
to $0 \leq 0$.  
For example, in network problems, the values of $\hat{a}_k(\cdot)$ and $\hat{b}_k(\cdot)$ 
often represent transmission rates, power expenditures, or newly accepted jobs, and the 
idle action is the one that accepts no new arrivals into the network and transmits no data over
any link of the network. For stock market problems, the idle action is often the 
action that neither buys nor sells any
shares of stock on the current slot.

\subsection{$T$-Slot Lookahead Policies} 

Rather than consider the optimum of the problem (\ref{eq:opt})-(\ref{eq:c4}) over the full time interval $t \in \{0, \ldots, t_{end}-1\}$, 
we consider the minimum cost that can be incurred over successive frames of size $T$, assuming that the 
time average constraints (\ref{eq:c1})-(\ref{eq:c3}) must be achieved over each frame.  
Specifically, 
let $T$ be a positive integer, representing a \emph{frame size}.   For a non-negative integer $r$, define
$F_r^*$ as the infimum value associated with the following problem (where $\bv{\gamma} \defequiv (\gamma_1, \ldots, \gamma_M)$): 
\begin{eqnarray}
\mbox{Minimize:} &  h_0  +  f(\bv{\gamma}) \label{eq:frame-lookahead-problem} \\
\mbox{Subject to:}  & h_l + g_l(\bv{\gamma}) \leq 0 \: \: \forall l \in \{1, \ldots, L\} \nonumber  \\
& \bv{\gamma} \in \script{X} \nonumber \\
 &\hspace{-.6in} \gamma_m = \frac{1}{T}\sum_{\tau=rT}^{(r+1)T-1} \hat{x}_m(\alpha(\tau), \omega(\tau)) \: \: \forall m\in\{1, \ldots, M\} \nonumber \\
 &\hspace{-.6in} h_l = \frac{1}{T}\sum_{\tau=rT}^{(r+1)T-1} \hat{y}_l(\alpha(\tau), \omega(\tau)) \: \: \forall l\in\{0, 1, \ldots, L\} \nonumber \\
 &\hspace{-.6in} \frac{1}{T}\sum_{\tau=rT}^{(r+1)T-1} [\hat{a}_k(\alpha(\tau), \omega(\tau)) - \hat{b}_k(\alpha(\tau), \omega(\tau))] \leq 0 \: \: \forall k \nonumber \\
 &\hspace{-1.1in} \alpha(\tau) \in \script{A}_{\omega(\tau)} \: \: \forall \tau \in \{rT, \ldots, (r+1)T-1\} \nonumber
\end{eqnarray} 
The value of $F_r^*$ represents the infimum of the cost metric 
that can be achieved over the frame, considering all policies that satisfy the constraints and 
that have perfect knowledge of the future $\omega(\tau)$ values over the frame. 
Our new goal is to design a non-anticipating control policy that is implemented over time $t_{end} = RT$ (for
some positive integer $R$), 
and that satisfies all constraints of the original problem while 
achieving a total cost that is close to (or smaller than) the value of: 
\begin{equation} \label{eq:val-T} 
 \frac{1}{R}\sum_{r=0}^{R-1} F_r^* 
 \end{equation} 
 
For $R=1$, it is clear that the value in (\ref{eq:val-T}) 
is the same as the optimal cost associated with the problem 
(\ref{eq:opt})  with $t_{end} = T$.  
For $R \geq 1$, it can be shown that the value in (\ref{eq:val-T}) 
is greater than or equal to the optimal cost associated with the problem 
(\ref{eq:opt}) for $t_{end} = RT$.  The reason that the problem  (\ref{eq:opt})-(\ref{eq:c4})  might have a strictly
smaller cost is that it only requires the time average constraints to be met over the full time
interval, rather than requiring them to be satisfied on each of the $R$ frames.  Nevertheless, when $T$ is large, 
it is not trivial to achieve the cost value of (\ref{eq:val-T}), as this cost is defined over policies that have $T$-slot
lookahead, whereas an actual policy does not have future lookahead capabilities.

\section{Solution to the General Problem} \label{section:general-solution}

First note that the problem (\ref{eq:opt})-(\ref{eq:c4}) is equivalent to the following problem, 
which introduces \emph{auxiliary variables} $\bv{\gamma}(t) = (\gamma_1(t), \ldots, \gamma_M(t))$ 
for $t \in \{0, \ldots, t_{end}-1\}$:
\begin{eqnarray}
\mbox{Minimize:} & \overline{y}_0 + \overline{f(\bv{\gamma})} \label{eq:new-opt} \\
\mbox{Subject to:} & \overline{y}_l + \overline{g_l(\bv{\gamma})} \leq 0 \: \: \forall l \in \{1, \ldots, L\} \label{eq:new-c1} \\
& \overline{a}_k \leq \overline{b}_k \: \: \forall k \in \{1, \ldots, K\} \label{eq:new-c2} \\
& \overline{\gamma}_m = \overline{x}_m \: \: \forall m \in \{1, \ldots, M\} \label{eq:new-c3} \\
& \alpha(t) \in \script{A}_{\omega(t)} \: \: \forall t \in \{0, \ldots, t_{end} -1\} \label{eq:new-c4} \\
&\bv{\gamma}(t) \in \script{X} \: \: \forall t \in \{0, \ldots, t_{end}-1\} \label{eq:new-c5} \\
& \hspace{-.4in} x_m^{min} \leq \gamma_m(t) \leq x_m^{max} \: \: \forall t \in \{0, \ldots, t_{end}-1\}  \label{eq:new-c6}
\end{eqnarray}
where $\overline{f(\bv{\gamma})}$ is defined: 
\[ \overline{f(\bv{\gamma})} \defequiv \frac{1}{t_{end}}\sum_{\tau=0}^{t_{end}-1} f(\bv{\gamma}(\tau)) \]
and where $\overline{g_l(\bv{\gamma})}$ is defined similarly. 

To see that the above problem (\ref{eq:new-opt})-(\ref{eq:new-c6}) 
is equivalent to the original problem (\ref{eq:opt})-(\ref{eq:c4}), note that any optimal solution of 
(\ref{eq:opt})-(\ref{eq:c4}) also satisfies the constraints (\ref{eq:new-c1})-(\ref{eq:new-c6}),  
with the same value of the cost metric (\ref{eq:new-opt}), 
provided that  we define $\gamma_m(t) = \overline{x}_m$ for all $t$, 
with $\overline{x}_m$ being the time average of $x_m(t)$ under the solution to the problem (\ref{eq:opt})-(\ref{eq:c4}).  
Thus, the minimum cost metric of the new problem (\ref{eq:new-opt})-(\ref{eq:new-c6}) 
is \emph{less than or equal} 
to that of (\ref{eq:opt})-(\ref{eq:c4}). 
On the other hand, by Jensen's inequality and convexity of $f(\bv{\gamma})$, $g_l(\bv{\gamma})$, 
we have for any solution of the new problem (\ref{eq:new-opt})-(\ref{eq:new-c6}): 
\begin{eqnarray*}
 \overline{f(\bv{\gamma})} &\geq& f(\overline{\bv{\gamma}}) =  f(\overline{\bv{x}}) \\
 \overline{g_l(\bv{\gamma})} &\geq& g_l(\overline{\bv{\gamma}}) = g_l(\overline{\bv{x}})
 \end{eqnarray*}
From this it easily follows that the minimum cost metric of the new problem (\ref{eq:new-opt})-(\ref{eq:new-c6}) 
is also \emph{greater than or equal} 
to that of the original problem (\ref{eq:opt})-(\ref{eq:c4}).

Such auxiliary variables are introduced in \cite{neely-fairness-infocom05} and \cite{now} to optimize
 functions of time averages, which is
very different from optimizing time averages of functions.  Note that if $f(\cdot) = g_l(\cdot) = 0$ for all $l$,
then $M=0$ and we do not need 
any auxiliary variables.\footnote{Also, in this case $M=0$, we do not need any 
queues of the type (\ref{eq:h-update}).}
Following the framework of \cite{neely-fairness-infocom05} \cite{now}, which solves problems similar to the above under ergodic assumptions on $\omega(t)$, in addition to the actual queues $Q_k(t)$ 
we define \emph{virtual queues} $Z_l(t)$ and $H_m(t)$ for each $l \in \{1, \ldots, L\}$
and $m \in \{1, \ldots, M\}$, with $Z_l(0) = H_m(0) = 0$ for all $l$ and $m$, and 
with update equation: 
\begin{eqnarray}
Z_l(t+1) &=& \max[Z_l(t) +  y_l(t) + g_l(\bv{\gamma}(t)), 0] \label{eq:z-update}   \\
H_m(t+1) &=&  H_m(t) + \gamma_m(t) - x_m(t)  \label{eq:h-update} 
\end{eqnarray}

Note that the queues $Q_k(t)$ and $Z_l(t)$ are non-negative for all $t$, while the queues
$H_m(t)$ can be possibly negative.  If the queues $Z_l(t)$, $Q_k(t)$,  $H_m(t)$ 
are close to zero at time $t_{end}$, then
the inequality constraints (\ref{eq:new-c1}), (\ref{eq:new-c2}), (\ref{eq:new-c3}), respectively,
are close to being satisfied, as specified by the following lemma.

\begin{lem} \label{lem:constraints} (Approximate Constraint Satisfaction) For any sequence $\{\omega(0), \omega(1), \omega(2), 
\ldots\}$, 
any designated end time $t_{end}>0$, any non-negative initial queue states $Z_l(0)$, $Q_k(0)$, any 
real valued initial queue states $H_m(0)$, and any sequence of control 
decisions $\alpha(0) \in \script{A}_{\omega(0)}$, $\alpha(1) \in \script{A}_{\omega(1)}$, $\alpha(2) \in \script{A}_{\omega(2)}$, etc., 
we have for all $l \in \{1, \ldots, L\}$, $k \in \{1, \ldots, K\}$:  
 \begin{eqnarray}
  \overline{a}_k &\leq& \overline{b}_k + \frac{Q_k(t_{end})- Q_k(0)}{t_{end}} \label{eq:foo2} \\
  \overline{y}_l + g_l(\overline{\bv{x}}) &\leq& \frac{Z_l(t_{end})-Z_l(0)}{t_{end}}  \nonumber\\
  && + \sum_{m=1}^M \frac{\beta_{l,m}|H_m(t_{end})-H_m(0)|}{t_{end}} \label{eq:truesat-g} \\
 &&\hspace{-.4in}\overline{\bv{x}} + \bv{\epsilon} \in \script{X} \label{eq:truesat-x} 
 \end{eqnarray}
 where $\bv{\epsilon} = (\epsilon_1, \ldots, \epsilon_M)$ is a vector with  components that satisfy: 
 \[ |\epsilon_m| = \frac{|H_m(t_{end}) - H_m(0)|}{t_{end}} \: \: \: \: \forall m \in \{1, \ldots, M\} \]
 where $\overline{\bv{x}} = (\overline{x}_1, \ldots, \overline{x}_M)$ represents
 a time average over the first $t_{end}$ slots, as defined in (\ref{eq:time-av-x}), 
 as do time averages $\overline{a}_k$, $\overline{b}_k$, $\overline{y}_l$. 
\end{lem} 

Note that the inequalities (\ref{eq:foo2})-(\ref{eq:truesat-x}) correspond to 
the desired constraints (\ref{eq:c1})-(\ref{eq:c3}), and show that these desired constraints are approximately 
satisfied if the values of $Q_k(t_{end})/t_{end}$, $Z_l(t_{end})/t_{end}$, $|H_m(t_{end})|/t_{end}$ are small. 
Recall that $\beta_{l,m}$ are  bounds on the partial derivatives of $g_l(\bv{x})$, and hence if for a particular $l$ we have
$g_l(\bv{x}) = 0$  for all $\bv{x}$, then $\beta_{l,m} = 0$ for all $m$, which tightens the bound in (\ref{eq:truesat-g}). 
This is useful for linear cost functions, as described in more detail in Section \ref{section:linear-costs}.

\begin{proof} (Lemma \ref{lem:constraints}) 
From the queueing update equation 
 (\ref{eq:q-update}) we have for any $k \in \{1, \ldots, K\}$ and  all $t\geq0$: 
 \[ Q_k(t+1) \geq Q_k(t) - b_k(t) + a_k(t) \]
 Summing over $\tau \in \{0, 1, \ldots, t-1\}$ for a given $t>0$, and dividing by $t$ gives: 
\begin{equation*} 
\frac{1}{t}\sum_{\tau=0}^{t-1} a_k(\tau) \leq \frac{1}{t}\sum_{\tau=0}^{t-1} b_k(\tau) + \frac{Q_k(t)-Q_k(0)}{t}
\end{equation*} 
Plugging $t = t_{end}$ into the above inequality proves (\ref{eq:foo2}). 

Similarly, the update equation (\ref{eq:h-update}) easily 
leads to the following for all $m \in \{1, \ldots, M\}$: 
\begin{equation*} 
  \left|\frac{1}{t}\sum_{\tau=0}^{t-1} \gamma_m(\tau) - \frac{1}{t}\sum_{\tau=0}^{t-1} x_m(\tau)\right| =  \frac{|H_m(t)- H_m(0)|}{t}
  \end{equation*}
Using $t=t_{end}$ in the above yields: 
\begin{equation} \label{eq:foo3} 
 |\overline{\gamma}_m - \overline{x}_m| = \frac{|H_m(t_{end}) - H_m(0)|}{t_{end}} 
 \end{equation} 
Note that $\bv{\gamma}(\tau) \in \script{X}$ for all $\tau \in \{0, \ldots, t_{end}-1\}$, and
hence (by convexity of $\script{X}$) the time average $\overline{\bv{\gamma}}$ is also in $\script{X}$.  Defining
$\bv{\epsilon} \defequiv \overline{\bv{\gamma}} - \overline{\bv{x}}$ thus ensures that $\overline{\bv{x}} + \bv{\epsilon} \in \script{X}$.
Noting that $\epsilon_m = \overline{\gamma}_m - \overline{x}_m$ and using (\ref{eq:foo3}) proves (\ref{eq:truesat-x}). 

Finally, 
from (\ref{eq:z-update}) we have for any $l \in \{1, \ldots, L\}$ and for 
all slots $t\geq0$: 
\[ Z_l(t+1) \geq Z_l(t) + y_l(t) + g_l(\bv{\gamma}(t)) \]
Summing over $\tau \in \{0, \ldots, t-1\}$ for any time $t>0$ yields: 
\[ Z_l(t) - Z_l(0) \geq \sum_{\tau=0}^{t-1} y_l(\tau)  + \sum_{\tau=0}^{t-1} g_l(\bv{\gamma}(\tau)) \]
Dividing by $t$ and rearranging terms proves that: 
\begin{equation} \label{eq:z-sat} 
 \frac{1}{t}\sum_{\tau=0}^{t-1} \left[y_l(\tau) + g_l(\bv{\gamma}(\tau))\right] \leq \frac{Z_l(t)- Z_l(0)}{t}  
 \end{equation} 
 Using $t=t_{end}$ in the above inequality, together with Jensen's inequality for the convex
 function $g_l(\cdot)$, yields: 
 \begin{equation} \label{eq:z-sat2} 
 \overline{y}_l + g_l(\overline{\bv{\gamma}}) \leq \frac{Z_l(t_{end})- Z_l(0)}{t_{end}}  
 \end{equation} 
 where $\overline{\bv{\gamma}}$ is the time average of $\bv{\gamma}(\tau)$ over the first $t_{end}$
 slots.  However, by (\ref{eq:cost-assumption-g}) we have: 
 \[   g_l(\overline{\bv{x}}) \leq g_l(\overline{\bv{\gamma}}) + \sum_{m=1}^M\beta_{l,m}|\overline{x}_m - \overline{\gamma}_m| \]
 Using this and (\ref{eq:foo3})  in (\ref{eq:z-sat2}) proves (\ref{eq:truesat-g}). 
\end{proof}

\subsection{Quadratic Lyapunov Functions and Sample Path Drift} 

Lemma \ref{lem:constraints} 
ensures the constraints (\ref{eq:new-c1})-(\ref{eq:new-c6}) 
(and hence the constraints  (\ref{eq:c1})-(\ref{eq:c4})) 
are approximately satisfied
if the final queue states of all queues are small relative to $t_{end}$.
Define $\bv{\Theta}(t)$ as a vector
of all current queue states: 
\begin{eqnarray*}
\bv{\Theta}(t) \defequiv [\bv{Z}(t), \bv{Q}(t), \bv{H}(t)] 
\end{eqnarray*}
where $\bv{Z}(t)$, $\bv{Q}(t)$, $\bv{H}(t)$ are vectors with entries $Z_l(t)$, $Q_k(t)$, $H_m(t)$, respectively. 
As a scalar 
measure of queue size,  define the following quadratic Lyapunov function, 
as in \cite{now}: 
\begin{equation} \label{eq:L}  
L(\bv{\Theta}(t)) \defequiv \frac{1}{2}\sum_{l=1}^L Z_l(t)^2 + \frac{1}{2}\sum_{k=1}^KQ_k(t)^2  + \frac{1}{2}\sum_{m=1}^MH_m(t)^2
\end{equation} 
For a given positive integer $T$, let $\Delta_T(t)$ 
represent the \emph{$T$-slot sample path Lyapunov drift} associated with particular
controls implemented over the interval $\{t, \ldots, t+T-1\}$ 
when the queues have state $\bv{\Theta}(t)$ at the beginning of the interval:\footnote{Note 
that the value of $\Delta_T(t)$ depends
on the queue state $\bv{\Theta}(t)$ at the start of the $T$-slot interval, the 
random events $\{\omega(t), \ldots, \omega(t+T-1)\}$, and the control actions
$\{\alpha(t), \ldots, \alpha(t+T-1)\}$ that are chosen over this interval.} 
\begin{equation} \label{eq:T-slot-drift-def}
 \Delta_T(t) \defequiv L(\bv{\Theta}(t+T)) - L(\bv{\Theta}(t)) 
 \end{equation} 
This notion of $T$-slot drift differs from that given in \cite{now} in that 
it does not involve an expectation.   It is difficult to control the $T$-slot drift, because
it depends on future (and hence unknown) $\omega(t)$ values.  Thus, following
the approach in \cite{now}, we design a control policy that, every slot $t$, 
observes the current $\omega(t)$ value
and the current queue states $\bv{\Theta}(t)$ 
and chooses a control action $\alpha(t) \in \script{A}_{\omega(t)}$ 
to minimize a weighted  sum of the 1-slot drift
and the current contribution to the cost metric (\ref{eq:new-opt}): 
\begin{eqnarray*}
\mbox{Minimize:} & \Delta_1(t) + V\hat{y}_0(\alpha(t), \omega(t)) + Vf(\bv{\gamma}(t))  \\
\mbox{Subject to:} & \mbox{Constraints (\ref{eq:new-c4})-(\ref{eq:new-c6})} 
\end{eqnarray*}
where $V \geq 0$ is a control parameter chosen in advance that affects a performance tradeoff. 
Rather than perform the exact minimization of the above problem, it suffices to minimize a bound. 
The following lemma bounds $\Delta_1(t)$. 

\begin{lem} \label{lem:lyap-drift} Under Assumption A1, the 1-slot drift $\Delta_1(t)$ satisfies: 
\begin{eqnarray*}
\Delta_1(t) \leq B + \sum_{l=1}^L Z_l(t)[y_l(t) +  g_l(\bv{\gamma}(t))] \\
+ \sum_{k=1}^KQ_k(t)[a_k(t) - b_k(t)] \\
 + \sum_{m=1}^MH_m(t)[\gamma_m(t) - x_m(t)] 
\end{eqnarray*}
where $B$ is a finite constant that satisfies  for all $t$:
\begin{eqnarray} 
B &\geq& \frac{1}{2}\sum_{l=1}^L(y_l(t) + g_l(\bv{\gamma}(t)))^2 + \frac{1}{2}\sum_{k=1}^K[b_k(t)^2 + a_k(t)^2] \nonumber \\
&& + \frac{1}{2}\sum_{m=1}^M (\gamma_m(t) - x_m(t))^2 \label{eq:B}
\end{eqnarray}
Such a finite constant $B$ exists by the boundedness assumptions A1, and a particular such $B$
is given  
in Appendix E. 
\end{lem} 
\begin{proof}
Squaring the $Z_l(t)$ update equation (\ref{eq:z-update}) and noting that $\max[x, 0]^2 \leq x^2$ yields: 
\begin{eqnarray}
 Z_l(t+1)^2 &\leq& Z_l(t)^2 + (y_l(t) + g_l(\bv{\gamma}(t)))^2  \nonumber \\
 && +2Z_l(t)[y_l(t) + g_l(\bv{\gamma}(t))]  \label{eq:one-slota}
 \end{eqnarray}
 Similarly, from (\ref{eq:q-update}) we have: 
 \begin{eqnarray}
 Q_k(t+1)^2 &\leq& (Q_k(t) - b_k(t))^2 + a_k(t)^2 \nonumber \\
 && + 2a_k(t)\max[Q_k(t) - b_k(t), 0] \nonumber \\
 &\leq& (Q_k(t) - b_k(t))^2 + a_k(t)^2  + 2a_k(t)Q_k(t)  \nonumber \\
 &=& Q_k(t)^2 + a_k(t)^2 + b_k(t)^2 \nonumber \\
 && +  2Q_k(t)[a_k(t) - b_k(t)] \label{eq:one-slotb} 
 \end{eqnarray}
 Finally, from (\ref{eq:h-update}) we have: 
 \begin{eqnarray} 
 H_m(t+1)^2 &=& H_m(t)^2 + (\gamma_m(t) - x_m(t))^2 \nonumber \\
 && + 2H_m(t)(\gamma_m(t) -x_m(t)) \label{eq:one-slotc} 
 \end{eqnarray} 
 Summing (\ref{eq:one-slota}),  (\ref{eq:one-slotb}), (\ref{eq:one-slotc}) and dividing by $2$ yields the result. 
\end{proof}

The above lemma shows that: 
\begin{eqnarray} 
\Delta_1(t) + V\hat{y}_0(\alpha(t), \omega(t)) + Vf(\bv{\gamma}(t)) \leq \nonumber \\
B + V\hat{y}_0(\alpha(t), \omega(t)) + Vf(\bv{\gamma}(t))  \nonumber \\
+ \sum_{l=1}^LZ_l(t)[\hat{y}_l(\alpha(t), \omega(t)) +g_l(\bv{\gamma}(t))] \nonumber \\
+  \sum_{k=1}^KQ_k(t)[\hat{a}_k(\alpha(t), \omega(t)) - \hat{b}_k(\alpha(t), \omega(t))] \nonumber \\
+ \sum_{m=1}^MH_m(t)[\gamma_m(t) - \hat{x}_m(\alpha(t), \omega(t))] \label{eq:big-1-drift} 
\end{eqnarray}

\subsection{The General Universal Scheduling Algorithm}  \label{section:general-alg} 

Our universal scheduling algorithm is designed to minimize the right hand side of (\ref{eq:big-1-drift}) 
every slot, as described  as follows: Every slot $t$, observe the random event 
$\omega(t)$ and the current queue backlogs $Z_l(t)$, $Q_k(t)$, $H_m(t)$ for all $l \in \{1, \ldots, L\}$, 
$k \in \{1, \ldots, K\}$, $m \in \{1, \ldots, M\}$, and perform the following. 
\begin{itemize} 
\item Choose $\bv{\gamma}(t) = (\gamma_1(t), \ldots, \gamma_M(t))$ to solve the following: 
\begin{eqnarray*}
\mbox{Minimize:} & Vf(\bv{\gamma}(t)) + \sum_{l=1}^LZ_l(t)g_l(\bv{\gamma}(t)) \\
& + \sum_{m=1}^MH_m(t)\gamma_m(t) \\
\mbox{Subject to:} & \bv{\gamma}(t) \in \script{X} \\
& x_m^{min} \leq \gamma_m(t) \leq x_m^{max} \: \: \forall m \in \{1, \ldots, M\} 
\end{eqnarray*}

\item Choose $\alpha(t) \in \script{A}_{\omega(t)}$ to minimize: 
\begin{eqnarray}
V\hat{y}_0(\alpha(t), \omega(t)) + \sum_{l=1}^LZ_l(t)\hat{y}_l(\alpha(t), \omega(t)) \nonumber  \\
- \sum_{m=1}^M H_m(t)\hat{x}_m(\alpha(t), \omega(t))\nonumber \\
 + \sum_{k=1}^KQ_k(t)[\hat{a}_k(\alpha(t), \omega(t)) - \hat{b}_k(\alpha(t), \omega(t))] \label{eq:expression} 
\end{eqnarray}

\item Update the actual queues $Q_k(t)$ and the virtual queues $Z_l(t)$, $H_m(t)$  for 
all $k \in \{1, \ldots, K\}$, $l \in \{1, \ldots, L\}$, $m \in \{1, \ldots, M\}$ via (\ref{eq:q-update}), (\ref{eq:z-update}), and (\ref{eq:h-update}).
\end{itemize} 

 Note that the above selection of $\bv{\gamma}(t)$ minimizes a convex 
function over a convex set, and decomposes into $M$ decoupled
convex optimizations of one variable in the case when cost functions $f(\bv{x})$ and $g_l(\bv{x})$ have  
the separable structure of (\ref{eq:separable-structure}) and when the set $\script{X}$ is equal to $\mathbb{R}^N$ or
a hypercube in $\mathbb{R}^N$.  
The optimization of $\alpha(t)$ in (\ref{eq:expression}) 
may be more complex and is possibly
a non-convex or combinatorial problem (depending on the $\script{A}_{\omega(t)}$ set and the 
$\hat{x}_m(\cdot)$, $\hat{y}_l(\cdot)$, $\hat{b}_k(\cdot)$, and $\hat{a}_k(\cdot)$ functions). However, it is simple when
the action space $\script{A}_{\omega(t)}$ contains only a finite (and small) number of options, in which case 
we can simply compare the functional (\ref{eq:expression}) for each option.
A key property of the above
algorithm is that it is \emph{non-anticipating} in that it acts only on the current $\omega(t)$ value, without
knowledge of future values. 

It can be shown that 
the expression (\ref{eq:expression}) has a well defined minimum value over the 
set $\script{A}_{\omega(t)}$ whenever Assumption A4 in Appendix B holds.  However, the
next theorem 
allows for \emph{approximate} minimization, 
where the choice of decision variables $\bv{\gamma}(t)$ and 
$\alpha(t)$ lead to a value that is off by an additive constant from achieving the minimum (or infimum) of the right hand
side in (\ref{eq:big-1-drift}). This is similar to the approximation results in \cite{now} and references therein, 
developed for ergodic problems. 
 Specifically, we define an algorithm to be \emph{$C$-approximate} if every slot it
makes decisions $\bv{\gamma}(t)$ and $\alpha(t)$ 
to satisfy the constraints (\ref{eq:new-c4})-(\ref{eq:new-c6}) and to yield either the infimum
of the expression on the right hand side of (\ref{eq:big-1-drift}) (as described in the algorithm above), or
to yield a value on the right hand side that differs from the infimum by at most an additive constant $C\geq 0$.

\begin{thm} \label{thm:1} Suppose Assumptions A1 and A2 hold. 
Consider any $C$-approximate algorithm. 
Let the random event sequence $\{\omega(0), \omega(1), \omega(2), \ldots \}$ be arbitrary.  Then:

(a) For any slot $t>0$ we have: 
\[ \sum_{l=1}^LZ_l(t)^2 + \sum_{k=1}^KQ_k(t)^2 + \sum_{m=1}^MH_m(t)^2 \leq tVC_0^2 + 2L(\bv{\Theta}(0))\]
where the constant $C_0$ is defined: 
\[ C_0 \defequiv \sqrt{2[(B + C)/V + (y_0^{max}-y_0^{min}) + (f^{max}  - f^{min})]} \]
where $B$ is the non-negative constant defined in (\ref{eq:B}).  In particular, all queues are bounded as follows: 
\[ Z_l(t), Q_k(t), |H_m(t)| \leq  \sqrt{tVC_0^2 + 2L(\bv{\Theta}(0))} \]
This bound becomes $C_0\sqrt{tV}$ if all queues are initially empty. 

(b)  If all queues are initially empty, then for any designated time $t_{end}>0$ we have:
\begin{eqnarray*}
& \overline{y}_l + g_l(\overline{\bv{x}}) \leq C_0\sqrt{\frac{V}{t_{end}}}\left[1 + \sum_{m=1}^M\beta_{l,m}\right] &  \forall l \in \{1, \ldots, L\} \\
& \overline{a}_k \leq \overline{b}_k + C_0\sqrt{\frac{V}{t_{end}}} &  \forall k \in \{1, \ldots, K\} \\
& \overline{\bv{x}} + \bv{\epsilon}(t) \in \script{X} & 
\end{eqnarray*} 
where $\bv{\epsilon}(t) = (\epsilon_1(t), \ldots, \epsilon_M(t))$ has entries that satisfy: 
\[ |\epsilon_m(t)| \leq C_0\sqrt{\frac{V}{t_{end}}} \] 
Thus, the constraints (\ref{eq:c4}) of the original problem are satisfied,  and the 
constraints (\ref{eq:c1})-(\ref{eq:c3}) of the original problem 
are approximately satisfied, where the error term in the approximation
decays with $t_{end}$ according to a constant multiple of $\sqrt{V/t_{end}}$.

(c) Consider any positive integer frame size $T$, any positive integer $R$, and define
$t_{end}  = RT$.   Then the value of the system cost
metric over $t_{end}$ slots satisfies: 
\begin{eqnarray}
 \overline{y}_0 + f(\overline{\bv{x}}) \leq \frac{1}{R}\sum_{r=0}^{R-1}F_r^* +  \frac{B+C}{V} 
+ \frac{D(T-1)}{V} \nonumber \\
 + \frac{L(\bv{\Theta}(0))}{VRT} +  \sum_{m=1}^M\frac{\nu_m|H_m(RT) - H_m(0)|}{RT}  \label{eq:utility-thm} 
\end{eqnarray}
where $\overline{y}_m$, $\overline{\bv{x}}$
are time averages over the first $t_{end}$ slots, and where 
$F_r^*$ is the optimal solution  to the problem (\ref{eq:frame-lookahead-problem})
and represents the optimal cost achieved by an idealized $T$-slot lookahead policy implemented
over the $r$th frame of size $T$. The constant $D$ is a finite constant that satisfies for all $t$
and all possible control actions that can be implemented on slot $t$:
\begin{eqnarray} 
D &\geq& \frac{1}{2}\sum_{l=1}^Lz_l^{diff}|y_l(t) + g_l(\bv{\gamma}(t))| \nonumber \\
&&  + \sum_{m=1}^Mh_m^{diff}|x_m(t) - \gamma_m(t)| \nonumber \\
&& + \frac{1}{2}\sum_{k=1}^K q_{k}^{diff}\max[b_k(t), a_k(t)]  \label{eq:D}
\end{eqnarray} 
where $z_l^{diff}$, $q_k^{diff}$, $h_m^{diff}$ represent the maximum change in queues
$Z_l(t)$, $Q_k(t)$, $H_m(t)$ over one slot, given by: 
\begin{eqnarray*}
z_l^{diff} &\defequiv& \max[|y_{l}^{max} + g_l^{max}|, |y_{l, min} + g_l^{min}|] \\
 q_k^{diff} &\defequiv& \max[b_k^{max}, a_k^{max}] \\
 h_m^{diff} &\defequiv& |x_m^{max} - x_m^{min}|
 \end{eqnarray*}
 Such a finite constant exists by the boundedness assumptions, and a value of $D$ that satisfies inequality 
 (\ref{eq:D}) is given in Appendix E. 

Finally, if initial queue backlogs are $0$, 
the final term in  (\ref{eq:utility-thm}) is bounded
by: 
\begin{eqnarray*}
 \sum_{m=1}^M\frac{\nu_m|H_m(RT) - H_m(0)|}{RT} \leq   \sum_{m=1}^M\nu_mC_0\sqrt{\frac{V}{RT}} 
\end{eqnarray*}
\end{thm} 
\begin{proof} 
See Appendix A.
\end{proof}

\subsection{Discussion of Theorem \ref{thm:1}} 
Consider the simple case when queues are initially 0, $g_l(\bv{x}) = 0$ for all $l$ (so that $\beta_{l,m} = 0$), 
and where we use a $C$-approximate algorithm for some
constant $C \geq 0$.  Fix a positive integer frame
size $T$. 
Theorem \ref{thm:1} can be interpreted as follows:   The algorithm implemented over $t_{end}=RT$ slots
ensures that the desired constraints (\ref{eq:c1})-(\ref{eq:c3}) are approximately met to within a ``fudge factor''
given by: 
\[ \mbox{Constraint fudge factor} = C_0\sqrt{\frac{V}{RT}} \]
where $C_0$ is the constant defined in part (a) of Theorem \ref{thm:1}. 
This fudge factor is made arbitrarily small 
when $RT$ is sufficiently larger than $V$.  Further, 
the achieved cost is either smaller than the cost $\frac{1}{R}\sum_{r=0}^{R-1} F_r^*$ associated with 
an ideal $T$-slot lookahead policy implemented over $R$ successive  frames of length $T$, or differs
from this value by an amount no more than a fudge factor that satisfies: 
\[ \mbox{Cost fudge factor} = \frac{C_1T}{V} + C_2 \sqrt{\frac{V}{RT}} \]
The value of $V$ can be chosen so that $C_1 T/V$ is arbitrarily small, in which case both the 
cost fudge factor and the constraint fudge factor are arbitrarily small provided that $R$ is sufficiently
large, that is, provided that we wait for a sufficiently large number of frames.

The constants $C_1$ and $C_2$ are given by: 
\begin{eqnarray*}
 C_1 &\defequiv& (B + C - D)/T + D  \\
C_2 &\defequiv& C_0 \sum_{m=1}^M \nu_m
\end{eqnarray*}
Note that $C_2=0$ in the case when $f(\bv{x})=0$ (so that $\nu_m = 0$ for all $m$).
Finally, note that the value of $T$ does not need to be chosen in order to implement the algorithm (we need
only choose a value of $V$), 
and hence the above bounds can be optimized over all positive integers $T$. 

The above tradeoffs described by $V$ and $R$ hold for general problems of the type (\ref{eq:opt})-(\ref{eq:c4}), 
and can be tightened for particular problems such as the network problem described in the next section, which 
provides queue bounds that do not grow with time. 
A similar strengthening due to constant queue bounds 
can be shown for the general problem in the case when Assumption A1 is strengthened
to a ``Slater-type'' condition, as described in Section \ref{section:slater}. 

\subsection{Linear Cost Functions} \label{section:linear-costs} 

The auxiliary variables are crucial for optimization of time varying systems
with non-linear cost functions $f(\bv{x})$, $g_l(\bv{x})$.
However, they are not needed when cost functions are linear (or affine).  For example, let $\bv{x}(t)$
be a vector of attributes as defined before, and consider the problem: 
\begin{eqnarray*}
\mbox{Minimize:} & h_0(\overline{\bv{x}}) \\
\mbox{Subject to:} & h_l(\overline{\bv{x}}) \leq 0 \: \: \:  \forall l \in \{1, \ldots, L\} \\
& \alpha(t) \in \script{A}_{\omega(t)}  \: \: \: \forall t \in \{0, \ldots, t_{end}-1\}
\end{eqnarray*}
where $h_l(\bv{x})$ are affine functions (i.e., linear plus a constant), so that $\overline{h_l(\bv{x})} = h_l(\overline{\bv{x}})$.
This can of course be treated using the framework of (\ref{eq:opt})-(\ref{eq:c4}) with $h_0(\bv{x}) = f(\bv{x})$, 
and $g_l(\bv{x}) = h_l(\bv{x})$, $y_l(\bv{x})=0$.   However, it can also be treated using (\ref{eq:opt})-(\ref{eq:c4}) with  
$f(\bv{x}) = g_l(\bv{x}) = 0$ for all $l \in \{1, \ldots, L\}$ and 
$y_l(\bv{x}) = h_l(\bv{x})$ (noting that $\overline{y}_l = h_l(\overline{x})$).  
This latter method is advantageous because it has $\nu_m = \beta_{l, m} = 0$ 
for all $l$ and $m$, which tightens the constraint inequalities (\ref{eq:truesat-g}) and the cost guarantee (\ref{eq:utility-thm}). 
Thus, it is useful to exploit linearity whenever possible. 

\subsection{Infinite Horizons and Ergodicity} \label{section:infinite-horizon}

Consider now the problem (\ref{eq:opt})-(\ref{eq:c4}) over an \emph{infinite horizon}, so that 
time averages $\overline{y}_l$, $\overline{\bv{x}}$, $\overline{a}_k$,
$\overline{b}_k$ represent limiting averages over the infinite horizon.   Suppose that Assumptions A1 and A2 hold, 
and that we use a $C$-approximate algorithm so that Theorem \ref{thm:1} applies.  
By taking a limit as $t_{end}\rightarrow\infty$, Theorem \ref{thm:1}b implies
that all required infinite horizon constraints are met.  Specifically, we have: 
\begin{eqnarray}
\limsup_{t\rightarrow\infty}\left[ \overline{y}_l(t) + g_l(\overline{\bv{x}}(t))\right] \leq 0 &  \forall l \in \{1, \ldots, L\} \label{eq:horizon1} \\
\limsup_{t\rightarrow\infty} \left[\overline{a}_k(t) - \overline{b}_k(t)\right] \leq 0 & \forall k \in \{1, \ldots, K\} \label{eq:horizon2} \\
\lim_{t_i\rightarrow\infty} \overline{\bv{x}}(t_i) \in \script{X} & \label{eq:horizon3} 
\end{eqnarray}
where $\overline{y}_l(t)$, $\overline{\bv{x}}(t)$, $\overline{a}_k(t)$, $\overline{b}_k(t)$ represent time
averages over the first $t$ slots, and where $\{t_i\}$ is any subsequence of times over which $\overline{\bv{x}}(t)$ converges.\footnote{Note
that $\overline{\bv{x}}(t)$ is an infinite sequence (with time index $t$) that takes values in a compact set, and so it has
a convergent subsequence.} 
Further, Theorem \ref{thm:1}c implies that the infinite horizon cost satisifes: 
\begin{eqnarray}
\hspace{-.3in}\lim_{R\rightarrow\infty} \left[\overline{y}_0(RT) + f\left(\overline{\bv{x}}(RT)\right)\right] &\leq& \lim_{R\rightarrow\infty}\frac{1}{R}\sum_{r=0}^{R-1} F_r^*
\nonumber \\
&& \hspace{-.1in}+ \frac{B + C}{V} + \frac{D(T-1)}{V} \label{eq:infinite-horizon}
\end{eqnarray}

Consider now the special case 
when the random events $\{\omega(0), \omega(1), \omega(2), \ldots\}$ evolve according to a 
general ergodic process with a well defined time average probability distribution.  In this case and 
under some mild assumptions, it can be shown that 
the optimal infinite horizon cost $f^*$ can be achieved over the class of stationary and randomized algorithms, 
that $F_r^*$ is close to $f^*$ for each $r$ whenever $T$ is sufficiently large, 
and that the term $\frac{1}{R}\sum_{r=0}^{R-1} F_r^*$ converges to $f^*$ plus an error term that is bounded by
$\delta(T)$, where $\delta(T)$ is a function that satisfies $\lim_{T\rightarrow\infty} \delta(T) = 0$. 
This is discussed in more detail in Appendix C.

\section{A Simple Internet Model} \label{section:internet}

Here we apply the universal scheduling framework to a simple flow based internet model, 
where we neglect the actual network queueing and develop a flow control policy that 
simply ensures the flow rate over any link is not more than the link capacity (similar to the flow
based models 
in \cite{kelly-shadowprice}\cite{low-flow-control}\cite{chiang-layering-decomposition}\cite{leonardi-SP-routing}\cite{chiang-delays-flownets}).  
Section \ref{section:network}  
treats a more extensive network model that explicitly  accounts for all queues. 

Suppose there are $N$ nodes and $L$ links, where each link $l \in \{1, \ldots, L\}$ 
has a possibly time-varying 
link capacity $C_l(t)$, for slotted time $t \in \{0, 1, 2, \ldots\}$. 
Suppose there are $M$ sessions, 
and let $A_m(t)$ represent the new arrivals to session $m$ on slot $t$.   Each session $m \in \{1, \ldots, M\}$
has a particular source node and a particular destination node. Assume link capacities
and newly arriving traffic are bounded so that: 
\begin{equation} \label{eq:internet-boundedness} 
 0 \leq C_l(t) \leq C_l^{max} \: \: \forall t \: \: , \: \: 0 \leq A_m(t) \leq A_m^{max} \: \: \forall t 
 \end{equation} 
for some finite constants $C_l^{max}$ and $A_m^{max}$. 
The random network event $\omega(t)$ is thus given by: 
\begin{equation} \label{eq:omega-internet} 
 \omega(t) \defequiv [(C_1(t), \ldots, C_L(t)); (A_1(t), \ldots, A_M(t))] 
 \end{equation} 
Recall that $\omega(t)$ is an arbitrary sequence with no probability model. 
The control action taken every slot is to first choose $x_m(t)$, the amount of type $m$ traffic admitted
into the network on slot $t$, according to: 
\[ 0 \leq x_m(t) \leq A_m(t) \]

Next, we must specify a \emph{path} for this data, from a collection of paths $\script{P}_m(t)$ associated
with path options of session $m$ on slot $t$ (possibly being the set of all possible paths in the network from the source
of session $m$ to its destination).\footnote{Strictly speaking, if there are time varying path choices then
we should augment $\omega(t)$ in (\ref{eq:omega-internet})  to include $\script{P}_m(t)$
for all $m\in\{1, \ldots, M\}$.}  Here, a path is defined in the usual sense, being a sequence of links starting
at the source, ending at the destination, and being such that the end node of each link is the start node of the next link. 
Let $1_{l,m}(t)$ be an indicator variable that is $1$ if the data $x_m(t)$ is selected to use a path that contains link $l$. 
The $(1_{l,m}(t))$ values completely specify the chosen paths for slot $t$, 
and hence the decision variable for slot $t$ is given by: 
\[ \alpha(t) \defequiv [(x_1(t), \ldots, x_M(t)); (1_{l,m}(t))|_{l \in \{1,\ldots, L\}, m \in \{1, \ldots, M\}}] \]

Let $\overline{\bv{x}} = (\overline{x}_1, \ldots, \overline{x}_M)$ be a vector 
of the infinite horizon time average admitted flow
rates, and let $\phi(\bv{x}) = \sum_{m=1}^M\phi_m(x_m)$ be a separable 
\emph{utility function}. Assume that each $\phi_m(x)$ is a continuous, 
concave, non-decreasing function in $x$, 
with maximum right derivative $\nu_m<\infty$.  Our goal is to maximize the throughput-utility 
$\phi(\overline{\bv{x}})$ subject to the constraints that the  time average flow over each link $l$ is
less than or equal to the time average capacity of that link. 
The infinite horizon utility optimization problem of interest
is thus: 
\begin{eqnarray*}
\mbox{Maximize:} & \sum_{m=1}^M \phi_m(\overline{x}_m) \\
\mbox{Subject to:} & \sum_{m=1}^M \overline{1_{l,m}x_m} \leq \overline{C}_l \: \: \: \forall l \in \{1, \ldots, L\}
\end{eqnarray*} 
where the time averages are defined: 
\begin{eqnarray*}
\overline{x}_m &\defequiv& \lim_{t\rightarrow\infty} \frac{1}{t}\sum_{\tau=0}^{t-1} x_m(\tau) \\
\overline{1_{l,m}x_m} &\defequiv& \lim_{t\rightarrow\infty} \frac{1}{t}\sum_{\tau=0}^{t-1} 1_{l,m}(\tau)x_m(\tau) \\
\overline{C}_l &\defequiv& \lim_{t\rightarrow\infty} \frac{1}{t}\sum_{\tau=0}^{t-1} C_l(\tau) 
\end{eqnarray*}
This is equivalent to minimizing the convex function $f(\overline{\bv{x}}) = -\phi(\overline{\bv{x}})$, and hence exactly fits our framework. 
There is no set constraint (\ref{eq:c3}), so that  $\script{X} = \mathbb{R}^M$.
As there are no actual queues $Q_k(t)$ in this model, we use only virtual queues $Z_l(t)$ and $H_m(t)$, defined by update
equations: 
\begin{eqnarray}
Z_l(t+1) = \max[Z_l(t) + \sum_{m=1}^M1_{l,m}(t)x_m(t) - C_l(t), 0] \label{eq:internet-z} \\
H_m(t+1) = H_m(t) + \gamma_m(t) - x_m(t) \label{eq:internet-h}
\end{eqnarray}
where $\gamma_m(t)$ are auxiliary variables for $m \in \{1, \ldots, M\}$.  
This is equivalent to the general framework with $y_l(t) = \sum_{m=1}^M1_{l,m}(t)x_m(t) - C_l(t)$ and $g_l(\cdot) = 0$
for $l \in \{1, \ldots, L\}$. 
Note that Assumption A1 holds by the boundedness assumptions (\ref{eq:internet-boundedness}), and Assumption A2 holds
because this system has an ``idle'' control action that admits no new data (so that $y_l(t) \leq 0$ for all $l$ under this
idle action). 
The general universal scheduling algorithm for this problem thus reduces to: 

\begin{itemize} 
\item (Auxiliary Variables) Every slot $t$, each session $m \in \{1, \ldots, M\}$ observes $H_m(t)$ 
and chooses $\gamma_m(t)$ as the solution 
to: 
\begin{eqnarray}
\mbox{Maximize:} & V\phi_m(\gamma_m(t)) - H_m(t)\gamma_m(t)  \label{eq:internet-aux1} \\
\mbox{Subject to:} & 0 \leq \gamma_m(t) \leq A_m^{max}  \label{eq:internet-aux2}
\end{eqnarray}
This is a simple maximization of a concave single variable function over an interval. 

\item (Routing and Flow Control) For each slot $t$ and each session $m \in \{1, \ldots, M\}$,
observe the new arrivals $A_m(t)$, the queue backlog $H_m(t)$, and the link 
queues $Z_l(t)$, and choose $x_m(t)$ and a path to maximize: 
\begin{eqnarray*}
\mbox{Maximize:} & x_m(t)H_m(t) - x_m(t)\sum_{l=1}^L 1_{l,m}(t)Z_l(t) \\
\mbox{Subject to:} & 0 \leq x_m(t) \leq A_m(t) \\
& \mbox{The path specified by $(1_{l,m}(t))$  is in $\script{P}_m(t)$}  
\end{eqnarray*}

This reduces to the following:  First find a shortest path from the source of session $m$
to the destination of session $m$, using link weights $Z_l(t)$ as link costs. If the total 
weight of the shortest path is less than or equal to $H_m(t)$, then choose $x_m(t) = A_m(t)$
and route all of this data over this single shortest path.  Else, there is too much congestion in the network,
and so we choose $x_m(t) = 0$ (thereby dropping
all data $A_m(t)$). 

\item (Virtual Queue Updates) Update the virtual queues according to (\ref{eq:internet-z}) and (\ref{eq:internet-h}). 
\end{itemize} 

The shortest path routing in this algorithm is similar to that given in \cite{leonardi-SP-routing}, which treats
a 
flow-based network stability problem 
in an ergodic setting under the assumption that arriving traffic is admissible (so that 
flow control is not used). 

Assume for simplicity that all queues are initially empty. Then 
for any frame size $T$ and any number of frames $R$, 
from Theorem \ref{thm:1}c, we know the utility of this algorithm satisfies: 
\begin{eqnarray}
\phi(\overline{\bv{x}}) &\geq& \frac{1}{R} \sum_{r=1}^{R-1}\Phi_r^* - \frac{B}{V} - \frac{D(T-1)}{V} \nonumber \\
&& - \sum_{m=1}^M\frac{\nu_m|H_m(RT)|}{RT} \label{eq:internet-util} 
\end{eqnarray}
where $\overline{\bv{x}}$ represents a time average over the first $RT$ slots, and $\Phi_r^*$ represents
the utility achieved by the $T$-slot lookahead policy implemented over slots $\{rT, \ldots, rT + T-1\}$. 
Here we are assuming we use an exact implementation of the algorithm (a $0$-approximation), 
so that $C=0$.  The constants $B$ and $D$ given in (\ref{eq:B}) and (\ref{eq:D}) are simplified 
for this context without queues $Q_k(t)$, and are provided in Appendix E.

Furthermore, the infinite horizon constraints are satisfied by (\ref{eq:horizon1}), 
and bounds on the virtual queue sizes for any time $t>0$ 
are also given in Theorem 1. 
However, with this particular structure we can obtain tighter bounds.  Indeed, from the update 
equation for $H_m(t)$ in (\ref{eq:internet-h}) and the auxiliary variable algorithm defined in (\ref{eq:internet-aux1})-(\ref{eq:internet-aux2}), 
it is easy
to see that: 
\begin{itemize} 
\item If $H_m(t) < 0$, then $\gamma_m(t) = A_m^{max}$ and hence $H_m(t)$ cannot decrease on the next slot. 

\item If $H_m(t) > V\nu_m$, then $\gamma_m(t) = 0$ and hence $H_m(t)$ cannot increase on the next slot. 
\end{itemize} 
It easily follows that: 
\begin{equation} \label{eq:easily-follows} 
 -A_m^{max}  \leq H_m(t) \leq V\nu_m + A_m^{max} \: \: \forall t \in \{0, 1, 2, \ldots\} 
 \end{equation} 
provided that this is true for $H_m(0)$ (which is indeed the case if $H_m(0) = 0$). 
Therefore, the final term in the utility guarantee (\ref{eq:internet-util}) is bounded by: 
\[ \sum_{m=1}^M \frac{\nu_m|H_m(RT)|}{RT} \leq \sum_{m=1}^M\frac{\nu_m(V\nu_m + A_m^{max})}{RT} \]
which goes to zero as the number of frames $R$ goes to infinity.  We further note that the utility 
guarantee (\ref{eq:internet-util}) can be modified to apply to \emph{any interval of $RT$ slots, starting at
any slot $t_0$}, provided that we modify the equation to account for possibly non-zero initial queue conditions
according to (\ref{eq:utility-thm}). 

Further, the fact that the queues $H_m(t)$ are deterministically
bounded allows one to deterministically bound the queue sizes $Z_l(t)$ as follows: For all $l \in \{1, \ldots, L\}$ we have: 
\begin{equation} \label{eq:z-net-bound} 
0 \leq Z_l(t)  \leq V\nu^{max} + (M+1)A^{max}  \: \: \: \: \forall t 
\end{equation} 
where $\nu^{max}$ and $A^{max}$ are defined as the maximum of all $\nu_m$ and $A_m^{max}$ values: 
\begin{eqnarray*}
\nu^{max} \defequiv \max_{m\in\{1, \ldots, M\}} \nu_m \: \: \: , \: \: \: 
A^{max} \defequiv \max_{m\in\{1, \ldots, M\}} A_m^{max} 
\end{eqnarray*}
The proof of this fact is simple:  If a link $l$ satisfies $Z_l(t) \leq V\nu^{max} + A^{max}$, then on the next slot
we have $Z_l(t+1) \leq V\nu^{max} + (M+1)A^{max}$, because the queue can increase by at most $MA^{max}$ on any
slot (see update equation (\ref{eq:internet-z})).  Else, if $Z_l(t) > V\nu^{max} + A^{max}$, then any path that uses this
link incurs a cost larger than $V\nu^{max} + A^{max}$, which is larger than $H_m(t)$ for any session $m$.  Thus, by the 
routing and flow control algorithm, no session will choose a path that uses this link 
on the current slot, and so $Z_l(t)$ cannot increase on the next
slot.    

\subsection{Delayed Feedback} 

We note that it may be difficult to use the exact queue values $Z_l(t)$ when solving for the shortest
path, as these values change every slot.  Hence, a practical implementation may use out-of-date
values $Z_l(t-\tau_{l,t})$
for some time delay $\tau_{l,t}$ that may depend on $l$ and $t$.  
Further, the virtual queue updates for $Z_l(t)$ in (\ref{eq:internet-z}) are most easily done at each link $l$, in which 
case 
the actual admitted data $x_m(t)$ for that link may not be known until some time delay, 
arriving as a process $x_m(t - \tau_{l,m.t})$. 
However, as the virtual queue size cannot change by more than a fixed
amount every slot, the queue value used differs from the ideal queue value by no more than an additive
constant that is proportional to the maximum time delay.  In this case, provided that the maximum time delay is bounded, 
we are simply using a $C$-approximation and the utility and queue bounds
are adjusted accordingly.  A more extensive treatment of delayed feedback for the case of 
networks without dynamic arrivals or channels is found in \cite{chiang-delays-flownets}, which 
uses a differential equation method. 

\subsection{Treating Wireless Networks with this Model} \label{section:internet-wireless} 

The above model can be applied equally to wireless networks.  However, an important extension in 
this case is to allow the link capacities $(C_1(t), \ldots, C_L(t))$ to be functions of a network resource
allocation decision (this is treated more extensively in Section \ref{section:network}).  This resource allocation
can be viewed as part of the network control action taken every slot. It is not difficult to show from the 
general solution in Section \ref{section:general-alg} that the optimal
resource allocation decision should observe $Z_l(t)$ values and 
choose capacities $C_l(t)$ to maximize the following weighted sum: 
\[ \sum_{l=1}^L C_l(t)Z_l(t) \]
Depending on the network model, this maximization can be difficult, and is generally prohibitively 
complex for wireless networks with interference.  Fortunately, it is easy to show that if the attempted max-weight
solution comes within a factor $\theta$ of the optimum max-weight decision every slot 
(for some value $\theta$ such
that $0 < \theta \leq 1$), then the same utility guarantees hold with $\Phi_r^*$ re-defined as the optimal $T$-slot
lookahead utility in a network with link capacities that are reduced by a factor $\theta$ from their actual values. 
This follows easily by noting that such a \emph{$\theta$-multiplicative-approximate} algorithm yields a right-hand-side
in the drift bound (\ref{eq:big-1-drift}) that is less than or equal to the right-hand-side associated with the optimal
max-weight decisions implemented on a network with $\theta$-reduced capacities (see also \cite{now} for 
a more detailed discussion of this for the case of i.i.d. $\omega(t)$ events).

\subsection{Limitations of this Model} 

The deterministic bound on $Z_l(t)$ in (\ref{eq:z-net-bound}) ensures that, 
over any interval of $K$ slots (for any positive integer $K$ and any initial slot $t_0$), 
the total data injected for use over link $l$ is no more than 
$V\nu^{max} + (M+1)A^{max}$ beyond the total capacity offered by the link: 
\begin{eqnarray*}
 \sum_{\tau=t_0}^{t_0 + K -1} \sum_{m=1}^M1_{l,m}(\tau)x_m(\tau) &\leq& \sum_{\tau=t_0}^{t_0+K-1} C_l(\tau) \\
 &&  + V\nu^{max} + (M+1)A^{max} 
 \end{eqnarray*}
 While this is a very strong deterministic bound that says no link is given more data than it can handle, it does not directly
 imply anything about the \emph{actual} network queues (other than the links are not overloaded).  The (unproven) 
 understanding is that, because the links are not overloaded, the actual network queues will be stable and all
 data can arrive to its destination with (hopefully small) delay. 

One might approximate average congestion or delay on a link as a convex function of the time average flow
rate over the link, as in \cite{pricing-congested-links}\cite{bertsekas-data-nets}\cite{chiang-delays-flownets}. 
This can be incorporated using the general framework of Section \ref{section:general-framework}, which allows for
optimization of time averages or convex functions of time averages.   However, we emphasize that this is only
an approximation and does not represent the actual network delay, or even a bound on delay. 
Indeed, while it is known that average queue congestion and delay is convex if a general 
stream of traffic is probabilistically split \cite{neely-convex-it}, this is
   not necessarily true (or relevant) for dynamically 
   controlled networks, particularly when the control depends on the queue backlogs and delays
   themselves.  Most problems involving optimization of actual network delay are difficult and unsolved. 
   Such problems 
   involve not only optimization of rate based utility functions, but engineering of the  
   Lagrange multipliers (which are related to queue backlogs) 
   associated with those utility functions.  

Finally, observe that the update equation for $Z_l(t)$ in (\ref{eq:internet-z}) 
can be interpreted as a queueing model where all admitted data on slot $t$ is placed immediately on 
all links $l$ of its path.  Similar models are used in 
\cite{low-flow-control}\cite{chiang-layering-decomposition}\cite{chiang-delays-flownets}\cite{lin-shroff-cdc04}. 
However, this is clearly an approximation, because data in an actual network will traverse
its path one link at a time. It is assumed that the actual network stamps all data with its intended path,
so that there is no dynamic re-routing mid-path.  Section \ref{section:network} treats an actual multi-hop queueing
network, 
and allows dynamic routing without pre-specified paths.

\section{Universal Network Scheduling} \label{section:network}

Consider a network with $N$ nodes that operates in slotted time.  
There are $M$ sessions, and we let $\bv{A}(t) = (A_1(t), \ldots, A_M(t))$ represent the vector of data that exogenously 
arrives to the transport layer for each session on slot $t$ (measured either
in integer units of \emph{packets} or real units of \emph{bits}).  
We assume that 
arrivals are bounded by  constants $A_{m}^{max}$, so that: 
\[ 0 \leq A_m(t) \leq A_{m}^{max} \: \: \forall t \]

Each session $m \in \{1, \ldots, M\}$ has a particular source node and 
destination node. 
 Data delivery takes place
by transmissions over possibly multi-hop paths.  We assume that a \emph{transport layer flow controller} observes
$A_m(t)$ every slot and decides how much of this data to add to the network layer at its source node, 
and how much to drop (flow control decisions are made to limit queue buffers and ensure the network is stable). 
Let $(x_m(t))|_{m=1}^M$ be the collection of \emph{flow control decision variables}  on 
slot $t$.  These decisions are made subject to the constraints: 
\begin{equation} \label{eq:flow-control-constraint} 
0 \leq x_m(t) \leq A_m(t) \: \: \forall m\in\{1, \ldots, M\}, \forall t
\end{equation} 

 All data that is intended for destination node
 $c \in \{1, \ldots, N\}$ is called \emph{commodity $c$ data}, regardless of its particular session. 
 For each $n \in \{1, \ldots, N\}$ and $c \in \{1, \ldots, N\}$, 
 let $\script{M}_n^{(c)}$ denote the set of all sessions $m \in \{1, \ldots, M\}$ that have source
 node $n$ and commodity $c$. 
 All data is queued according to its
commodity, and we define $Q_n^{(c)}(t)$ as
the amount of commodity $c$ data in node $n$ on slot $t$. We assume that $Q_n^{(n)}(t) = 0$ for all $t$, 
as data that reaches its destination is removed from the network.  
Let $\bv{Q}(t)$ denote the matrix of current
queue backlogs for all nodes and commodities.  

The queue backlogs change from slot to slot as follows: 
\[ Q_n^{(c)}(t+1) = Q_n^{(c)}(t) - \sum_{j=1}^N \tilde{\mu}_{nj}^{(c)}(t) + \sum_{i=1}^N \tilde{\mu}_{in}^{(c)}(t) + \sum_{m\in\script{M}_n^{(c)}}x_m(t) \]
where $\tilde{\mu}_{ij}^{(c)}(t)$ denotes the actual amount of commodity $c$ data transmitted from node 
$i$ to node $j$ (i.e., over link $(i,j)$) on slot $t$. It is useful to define \emph{transmission decision variables}
$\mu_{ij}^{(c)}(t)$ as the bit rate \emph{offered} by  link $(i,j)$ to commodity $c$ data, where this full amount is
used if there is that much commodity $c$ data available at node $i$, so that: 
\[ \tilde{\mu}_{ij}^{(c)}(t) \leq  \mu_{ij}^{(c)}(t) \: \: \forall i,j, c \in \{1, \ldots, N\}, \forall t \]
For simplicity, we assume that if there is not enough data to send at the offered rate, then \emph{null data} is sent, so
that:\footnote{All results hold exactly as stated if this null data is not sent, because the drift
bound in Lemma \ref{lem:lyap-drift} holds exactly wen the update equation (\ref{eq:q-net}) is replaced
by an inequality $\leq$, see \cite{now}.}  
\begin{eqnarray} 
Q_n^{(c)}(t+1) &=& \max[Q_n^{(c)}(t) - \sum_{j=1}^N\mu_{nj}^{(c)}(t), 0] \nonumber \\
&& + \sum_{i=1}^N\mu_{in}^{(c)}(t) + \sum_{m \in \script{M}_n^{(c)}} x_m(t) \label{eq:q-net} 
\end{eqnarray}

This satisfies (\ref{eq:q-update}) if we relate index $k$ (for $Q_k(t)$ in (\ref{eq:q-update})) to index $(n,c)$ (for 
$Q_n^{(c)}(t)$ in (\ref{eq:q-net})), and if we define: 
\begin{eqnarray*}
b_n^{(c)}(t) &\defequiv& \sum_{j=1}^N \mu_{nj}^{(c)}(t) \\
a_n^{(c)}(t) &\defequiv& \sum_{i=1}^N \mu_{in}^{(c)}(t) + \sum_{m \in \script{M}_n^{(c)}} x_m(t)
\end{eqnarray*}

\subsection{Transmission Variables} 
Let $S(t)$ represent the \emph{topology state} of the network on slot $t$,  observed on each slot $t$
as in \cite{now}.  
The value of $S(t)$ is an abstract and possibly multi-dimensional quantity that describes the current
link conditions between all nodes under the current slot.  The collection of all transmission rates that can be 
offered over each link $(i,j)$ of the network is given by a general transmission rate function $\bv{C}(I(t), S(t))$:\footnote{It is 
worth noting now that for networks with orthogonal channels, 
our ``max-weight'' transmission algorithm (to be defined in the next subsection) 
decouples to allow nodes to make transmission decisions 
based only on those components of the current topology state $S(t)$ that relate to their
own local channels.
Of course, for wireless interference networks, all channels are coupled,  although distributed 
approximations of max-weight transmission exist in this case \cite{now}.} 
\[ \bv{C}(I(t), S(t)) = (C_{ij}(I(t), S(t)))_{i,j\in\{1, \ldots, N\}, i\neq j} \]
where $I(t)$ is a general network-wide resource allocation decision (such as link scheduling, bandwidth 
selection, modulation, etc.) and takes values in some abstract set $\script{I}_{S(t)}$ that possibly depends
on the current $S(t)$. 
We assume that the transmission rate function $C_{ij}(I(t), S(t))$ is non-negative and bounded by a finite
constant $\mu_{ij}^{max}$ for all $(i,j)$, $I(t)$, and $S(t)$. 

Every slot the network controller observes the current $S(t)$ and 
makes a resource allocation decision $I(t) \in \script{I}_{S(t)}$.  The controller then chooses
$\mu_{ij}^{(c)}(t)$ variables subject to the following constraints: 
\begin{eqnarray}
\mu_{ij}^{(c)}(t) \geq 0 &  \forall i,j,c \in\{1, \ldots, N\} \label{eq:mu1} \\
\mu_{ii}^{(c)}(t) = \mu_{ij}^{(i)}(t) = 0 & \forall i, j, c \in \{1, \ldots, N\} \label{eq:mu2}  \\
\sum_{c=1}^N \mu_{ij}^{(c)} \leq C_{ij}(I(t), S(t)) &  \forall i,j \in \{1, \ldots, N\} \label{eq:mu3} 
\end{eqnarray}

\subsection{The Utility Optimization Problem} 

This problem fits the general model of Section \ref{section:general-framework} 
by defining the random event $\omega(t)$ as
follows: 
\[ \omega(t) \defequiv [\bv{A}(t); S(t)] \]
That is, the random event $\omega(t)$ is the collection of all new arrivals together with the current topology state. 
The control action $\alpha(t)$ is defined by: 
\[ \alpha(t) \defequiv [I(t); (\mu_{ij}^{(c)}(t))|_{i,j,c \in \{1, \ldots, N\}}; (x_m(t))|_{m=1}^M] \]
representing the resource allocation, transmission, and flow control decisions.  The action space $\script{A}_{\omega(t)}$ 
is defined by the set of all $I(t) \in \script{I}_{S(t)}$, all $(\mu_{ij}^{(c)}(t))$ that satisfy (\ref{eq:mu1})-(\ref{eq:mu3}), 
and all $(x_m(t))$ that satisfy (\ref{eq:flow-control-constraint}). 

Define $\overline{x}_m$ as the time average of $x_m(t)$ over the first
$t_{end}$ slots (as in (\ref{eq:time-av-x})), and 
define $\overline{\bv{x}}$ as the vector of these time averages. 
Our objective is to 
solve the following problem: 
\begin{eqnarray}
\mbox{Maximize:} &  \phi(\overline{\bv{x}}) \label{eq:net-opt1} \\
\mbox{Subject to:} &  \alpha(t) \in \script{A}_{\omega(t)} \: \: \forall t \in \{0, \ldots, t_{end}-1\} \label{eq:netc1}   \\
& \sum_{m\in\script{M}_n^{(c)}} \overline{x}_m + \sum_{i=1}^N \overline{\mu}_{in}^{(c)} \leq
 \sum_{j=1}^N \overline{\mu}_{nj}^{(c)}  \nonumber \\
 & \hspace{+1.3in} \forall n,c \in \{1, \ldots, N\}  \label{eq:netc2} 
\end{eqnarray}
where $\phi(\overline{\bv{x}})$ is a continuous, concave, and entrywise non-decreasing  
utility function of the form: 
\[ \phi(\bv{x}) \defequiv \sum_{m=1}^M \phi_m(x_m) \]
Define $\nu_m$ as the right partial derivative of $\phi_m(x)$ at $x=0$, and assume
$0 \leq \nu_m < \infty$ for all $m$. 
Thus, this problem fits exactly into the general framework, satisfying Assumption A1 by our 
boundedness assumptions, and satisfying Assumption A2 by the ``idle'' control action that
admits no new data and transmits over no links. 

\subsection{The Universal Network Scheduling Algorithm} 

To apply the general solution, note that the constraints (\ref{eq:netc2}) are upheld by stabilizing the actual
queues $Q_n^{(c)}(t)$ with updates (\ref{eq:q-net}).  Because we have not specified any additional constraints, 
there are no $Z_l(t)$ queues as used in the general framework.  However, we have auxiliary variables
$\gamma_m(t)$ for each $m \in \{1, \ldots, M\}$, with virtual queue update: 
\begin{equation} \label{eq:h-net-update} 
H_m(t+1) = H_m(t) + \gamma_m(t) - x_m(t) 
 \end{equation} 
The algorithm is thus:  
\begin{itemize} 
\item (Auxiliary Variables) For each slot $t$, each session $m \in \{1, \ldots, M\}$ observes 
the current virtual queue $H_m(t)$, and chooses auxiliary variable
$\gamma_m(t)$ as the solution to: 
\begin{eqnarray}
\mbox{Maximize:} & V\phi_m(\gamma_m(t)) - H_m(t)\gamma_m(t) \label{eq:aux1} \\
\mbox{Subject to:} & 0 \leq \gamma_m(t) \leq A_{m}^{max} \nonumber  
\end{eqnarray}
This is a maximization of a concave single variable function over an interval, the same
as in the internet algorithm of Section \ref{section:internet}. 

\item (Flow Control) For each slot $t$, each session $m$ observes $A_m(t)$ and
the queue
values $H_m(t)$,  $Q_{n_m}^{(c_m)}(t)$ (where $n_m$ denotes the source
node of session $m$, and $c_m$ represents its destination).  Note that these queues
are all local to the source node of the session, and hence can be observed easily.  
It then chooses $x_m(t)$  to solve: 
\begin{eqnarray}
\mbox{Maximize:} & H_m(t)x_m(t) - Q_{n_m}^{(c_m)}(t)x_m(t)  \label{eq:flow-control} \\
\mbox{Subject to:} & 0 \leq x_m(t) \leq A_m(t)  \nonumber 
\end{eqnarray}
This reduces to the ``bang-bang'' flow control decision of 
choosing $x_m(t) = A_m(t)$ if $Q_{n_m}^{(c_m)}(t) \leq H_m(t)$, and $x_m(t) = 0$
otherwise. 

\item (Resource Allocation and Transmission) For each slot $t$, the network controller observes
queue backlogs $\{Q_n^{(c)}(t)\}$ and the topology state $S(t)$ and chooses $I(t) \in \script{I}_{S(t)}$ and 
$\{\mu_{ij}^{(c)}(t)\}$ to solve: 
\begin{eqnarray}
\mbox{Max:} & \hspace{-.1in}\sum_{n,c} Q_n^{(c)}(t)[\sum_{j=1}^N\mu_{nj}^{(c)}(t) - \sum_{i=1}^N\mu_{in}^{(c)}(t)   ] \label{eq:max-weight}\\
\mbox{S.t.:} &  I(t) \in \script{I}_{S(t)} \mbox{ and (\ref{eq:mu1})-(\ref{eq:mu3})} \nonumber
\end{eqnarray}

\item (Queue Updates) Update the virtual queues $H_m(t)$ according to (\ref{eq:h-net-update}) and
the acutal queues $Q_n^{(c)}(t)$ according to  (\ref{eq:q-net}). 
\end{itemize} 

The exact decisions required to implement the resource allocation and transmission component are described
in Subsection \ref{subsection:transmission} below.  Before covering this, we state the performance of the algorithm 
under a general $C$-approximate implementation of the above algorithm.  For simplicity, assume
all queues are initially zero.  By Theorem \ref{thm:1}c we have
for any integers $R>0$, $T>0$: 
\begin{eqnarray}
\phi(\overline{\bv{x}}) &\geq& \frac{1}{R}\sum_{r=0}^{R-1} \Phi_r^* - \frac{B+C}{V} - \frac{D(T-1)}{V} \nonumber \\ 
&& - \sum_{m=1}^M\frac{\nu_m|H_m(RT)|}{RT}   \label{eq:util-net} 
\end{eqnarray}

While Theorem \ref{thm:1} also provides a bound on the final term, and bounds on all queue sizes, we can again
provide tighter constant queue 
bounds by taking advantage of the flow control structure of the problem. Indeed, by the same argument that proves (\ref{eq:easily-follows}) 
in Section \ref{section:internet}, we have that $H_m(t)$ cannot decrease if it is already negative, and cannot increase
if it is beyond $V\nu_m$, so that for all $m \in \{1, \ldots, M\}$ we have: 
\begin{equation} \label{eq:h-net-bound} 
 -A_m^{max} \leq H_m(t) \leq V\nu_m + A_m^{max} \: \: \: \forall t
 \end{equation} 
provided that these bounds are true for $H_m(0)$ (which is indeed the case if $H_m(0) = 0$). 
Therefore, $|H_m(t)| \leq V\nu_m + A_m^{max}$, and
the utility bound (\ref{eq:util-net}) becomes: 
\begin{eqnarray}
\phi(\overline{\bv{x}}) &\geq& \frac{1}{R}\sum_{r=0}^{R-1} \Phi_r^* - \frac{B+C}{V} - \frac{D(T-1)}{V} \nonumber \\ 
&& - \sum_{m=1}^M\frac{\nu_m(V\nu_m + A_m^{max})}{RT}   \label{eq:util-net2} 
\end{eqnarray}
The values of $B$ and $D$ in (\ref{eq:util-net2}) for this context are given in Appendix E.  Note that this yields
a ``utility fudge factor'' of the form as indicated in the introduction of this paper: 
\[ \mbox{utility fudge factor} = \frac{B_1T}{V} + \frac{B_2V}{RT} \]
where: 
\begin{eqnarray*}
B_1 \defequiv (B + C - D)/T + D \: \: , \: \: 
B_2 \defequiv \sum_{m=1}^M\nu_m(\nu_m + A_m^{max}/V)
\end{eqnarray*}

The value of $C$ used in the above bound is equal to $0$ if we use a $0$-approximation, being an exact implementation
of the above algorithm.  In the next subsections we purposefully engineer a $C$-approximation for a nonzero but constant 
$C$ so that we can additionally provide deterministic bounds on all \emph{actual} queues $Q_n^{(c)}(t)$.

\subsection{Resource Allocation and Transmission} \label{subsection:transmission} 

By switching the sums in (\ref{eq:max-weight}), 
it is easy to show that the resource allocation and transmission maximization
reduces to the following generalized ``max-weight'' algorithm 
(see \cite{now}):  Every slot $t$, choose $I(t) \in \script{I}_{S(t)}$ to maximize: 
\begin{eqnarray*}
 \sum_{i=1}^N \sum_{j=1}^N C_{ij}(I(t), S(t))W_{ij}(t)
\end{eqnarray*}
where $W_{ij}(t)$ are weights defined by: 
\begin{eqnarray*}
W_{ij}(t) \defequiv \max_{c \in \{1, \ldots, N\}} \max[W_{ij}^{(c)}(t), 0]
\end{eqnarray*}
where $W_{ij}^{(c)}(t)$ are differential backlogs: 
\[ W_{ij}^{(c)}(t) \defequiv Q_i^{(c)}(t) - Q_j^{(c)}(t) \]
The transmission decision variables are then given by: 
\[ \mu_{ij}^{(c)}(t) = \left\{ \begin{array}{ll}
                             C_{ij}(I(t), S(t)) &\mbox{ if $c = c^*_{ij}(t)$ and $W_{ij}^{(c)}(t) \geq 0$} \\
                             0  & \mbox{ otherwise} 
                            \end{array}
                                 \right.  \]
where $c_{ij}^*(t)$ is defined as the commodity $c \in \{1, \ldots, N\}$
that maximizes the differential backlog $W_{ij}^{(c)}(t)$ (breaking ties arbitrarily). 

\subsection{A $C$-Approximate Transmission Algorithm} \label{subsection:c-approximate-transmission} 

Rather than implement the exact transmission algorithm in the above subsection, 
we present here a useful $C$-approximation that yields bounded queues (see also 
\cite{neely-mesh} \cite{neely-energy-it}). Define $\hat{W}_{ij}^{(c)}(t)$ as follows: 
\begin{eqnarray} 
 & \hspace{-2.8in} \hat{W}_{ij}^{(c)}(t) \defequiv  \nonumber \\
& \left\{ \begin{array}{ll}
                         W_{ij}^{(c)}(t) + \theta_i^{(c)} - \theta_j^{(c)} & \mbox{if $Q_j^{(c)}(t) \leq Q^{max} - \beta_j$}    \\
                             -1  & \mbox{ otherwise} 
                            \end{array}
                                 \right.  \label{eq:new-weights} 
\end{eqnarray} 
where for each $n \in \{1, \ldots, N\}$, 
$\beta_n$ is defined as the largest amount of any commodity that can enter node $n$, considering both exogenous
and endogenous arrivals (this is finite by the boundedness assumptions on transmission rates and new arrivals), 
and where $Q^{max}$ is defined: 
\begin{eqnarray}
Q^{max} \defequiv V\nu^{max} + A^{max} + \beta^{max} \label{eq:qmax} 
\end{eqnarray} 
where $\nu^{max}$, $A^{max}$, and $\beta^{max}$ are given by: 
\begin{eqnarray*}
\nu^{max} &\defequiv& \max_{m\in\{1, \ldots, M\}} \nu_m \\
A^{max} &\defequiv& \max_{m\in\{1, \ldots, M\}} A_m^{max} \\
\beta^{max} &\defequiv& \max_{n \in \{1, \ldots, N\}} \beta_n
\end{eqnarray*}
Finally, the values $\theta_i^{(c)}$ are any non-negative weights that represent some type of estimate
of the distance from node $i$ to destination $c$ (possibly being zero if there is no such estimate available). 
Such weights are known to experimentally improve delay by 
biasing routing decisions to move in directions closer
to the destination (see \cite{neely-power-network-jsac}\cite{now}\cite{neely-mesh}). 
Then define $\hat{W}_{ij}(t)$ as: 
\[ \hat{W}_{ij}(t) \defequiv \max_{c \in \{1, \ldots, N\}} \max[\hat{W}_{ij}^{(c)}(t), 0] \]
and choose $I(t) \in \script{I}_{\omega(t)}$ to maximize: 
\begin{equation} \label{eq:max-weight-enhanced} 
 \sum_{i=1}^N\sum_{j=1}^N C_{ij}(I(t), S(t)) \hat{W}_{ij}(t)  
 \end{equation} 
and choose transmission variables: 
\begin{equation}  \mu_{ij}^{(c)}(t) =  \left\{ \begin{array}{ll}
                             C_{ij}(I(t), S(t)) &\mbox{ if $c = \hat{c}^*_{ij}(t)$ and $\hat{W}_{ij}^{(c)}(t) \geq 0$} \\
                             0  & \mbox{ otherwise} 
                            \end{array}
                                 \right. \label{eq:mu-ij-enhanced}  
\end{equation} 
where $\hat{c}^*(t)$ is the commodity $c \in \{1, \ldots, N\}$ 
that maximizes $\hat{W}_{ij}^{(c)}(t)$.

\subsection{Bounded Queues}

\begin{lem} (Bounded $Q_n^{(c)}(t)$) Suppose auxiliary variables and flow control decisions are made according 
to (\ref{eq:aux1}) and (\ref{eq:flow-control}), with update equations (\ref{eq:h-net-update}) and (\ref{eq:q-net}). 
Suppose that $I(t) \in \script{I}_{S(t)}$ is chosen in some arbitrary
manner every slot $t$ (not necessarily according to (\ref{eq:max-weight-enhanced})), but that
transmission decisions are made according to (\ref{eq:mu-ij-enhanced}) with respect to the particular
$I(t)$ chosen.  Then for all $t$ we have: 
\begin{eqnarray}
Q_n^{(c)}(t) \leq Q^{max}  = V\nu^{max} + A^{max} + \beta^{max} \label{eq:det-q-bound}
\end{eqnarray}
provided that this inequality holds at $t=0$. 
\end{lem}

\begin{proof} 
Suppose that $Q_n^{(c)}(t) \leq Q^{max}$ for all $n,c$ for a particular slot $t$ (this is true
by assumption on slot $t=0$).  We prove it also holds for slot $t+1$.  First suppose that  
$Q_n^{(c)}(t) \leq Q^{max} - \beta_n$.  Then, because $\beta_n$ is the largest amount of new arrivals to queue
$Q_n^{(c)}(t)$ over one slot (considering both endogenous and exogenous arrivals), it must be that
$Q_n^{(c)}(t+1) \leq Q^{max}$.  

Consider now the opposite case when  
$Q^{max} - \beta_n < Q_n^{(c)}(t) \leq Q^{max}$. Then from (\ref{eq:new-weights}) 
we see that $\hat{W}_{in}^{(c)} = -1$ for all links $(i,n)$ over which new commodity $c$ data
could be transmitted to node $n$ from other nodes. Thus, by (\ref{eq:mu-ij-enhanced}) we see that
no commodity $c$ data will be transmitted to node $n$ from any other node.  Further, 
We have: 
\[ Q_n^{(c)}(t) > Q^{max} - \beta_n = V\nu^{max} + A^{max}  \geq H_m(t) \]
for all $m \in \{1, \ldots, M\}$, 
where the first equality follows by definition of $Q^{max}$ in 
(\ref{eq:qmax}) and the 
final inequality follows by (\ref{eq:h-net-bound}).  It follows by the flow control
decision (\ref{eq:flow-control}) that $x_m(t) = 0$ for all sessions $m$ that might deliver 
new data to queue $Q_{n}^{(c)}(t)$.  Thus, 
no new commodity $c$
data (exogenous or endogenous)
arrives to node $n$ on slot $t$, and $Q_n^{(c)}(t+1) \leq Q_{n}^{(c)}(t) \leq Q^{max}$. 
\end{proof} 

The queue bound (\ref{eq:det-q-bound}) in the above lemma provides the strong deterministic
guarantee that all queues are bounded by a constant that grows linearly with the $V$ parameter. 
Thus, while increasing $V$ can improve the terms $(B+C)/V$ and $D(T-1)/V$ in the utility guarantee
(\ref{eq:util-net2}), a tradeoff is in the linear growth with $V$ in queue congestion (\ref{eq:det-q-bound}), 
as well as the increase in the number of frames $R$ required for the final term in the utility bound (\ref{eq:util-net2})
to decay to near-zero. 

While it is intuitive that the above algorithm produces a $C$-approximation for some constant value $C$, we complete
the analysis below by formally showing this.  Additionally, we note that a $\theta$-multiplicative-approximate 
solution to (\ref{eq:max-weight-enhanced})
leads to utility guarantees where $\Phi_r^*$ is re-defined as a $T$-slot lookahead utility on a network where 
the $\bv{C}(I(t), S(t))$ function is replaced by $\theta\bv{C}(I(t), S(t))$, which holds for the same reason as described in 
Section \ref{section:internet-wireless}. 

\subsection{Computing the $C$ value} 

\begin{lem} Using the modified weights $\hat{W}_{ij}^{(c)}(t)$ in (\ref{eq:new-weights}) results in a $C$-approximation
of the max-weight resource allocation and transmission scheduling problem (\ref{eq:max-weight}), with: 
\begin{equation} \label{eq:C}
C \defequiv 2C_{sum}[\beta^{max} + \theta_{diff}]
\end{equation} 
where $C_{sum}$ is the largest possible sum of transmission rates $\sum_{ij} C_{ij}(I(t), S(t))$, summed over all links
and considering all possible $S(t)$ states and $I(t)$ decisions (being finite by the boundedness assumptions on all links), 
and $\theta_{diff}$ is the maximum difference in $\theta_i^{(c)}$ and $\theta_j^{(c)}$, maximized over all node pairs $(i,j)$
and all commodities $c$. 
\end{lem} 
\begin{proof} 
Because all queues $Q_n^{(c)}(t)$ are upper bounded by $Q^{max}$, if $Q_j^{(c)}(t) > Q^{max} - \beta_j$, then
$\max[W_{ij}^{(c)}(t), 0] = \max[Q_i^{(c)}(t) - Q_j^{(c)}(t), 0] \leq \beta_j$.  It follows that: 
\[ |\max[W_{ij}^{(c)}(t), 0] - \max[\hat{W}_{ij}^{(c)}(t), 0] | \leq \beta^{max} + \theta_{diff} \]
It follows that: 
\[ |W_{ij}(t) - \hat{W}_{ij}(t)| \leq \beta^{max} + \theta_{diff} \]
Therefore: 
\begin{eqnarray*}
|\sum_{i=1}^N\sum_{j=1}^N C_{ij}(I(t), S(t)) [W_{ij}(t) - \hat{W}_{ij}(t)]|  \\
\leq \sum_{i=1}^N \sum_{j=1}^N C_{ij}(I(t), S(t)) [\beta^{max} + \theta_{diff}] \\
\leq  C_{sum}[\beta^{max} + \theta_{diff}]  = C/2
\end{eqnarray*}
where $C_{sum}$ is the maximum sum rate over all links on any slot. 
Now let $I^*(t)$ be the maximum of $\sum_{ij} C_{ij}(I(t), S(t))W_{ij}(t)$ over
$\script{I}_{S(t)}$, and let $\hat{I}(t)$ be the maximum of 
$\sum_{ij}C_{ij}(I(t), S(t))\hat{W}_{ij}(t)$. Then: 
\begin{eqnarray*}
&& \hspace{-1in}\sum_{ij} C_{ij}(\hat{I}(t), S(t)) W_{ij}(t) \\
&\geq& 
\sum_{ij}C_{ij}(\hat{I}(t), S(t)) \hat{W}_{ij}(t) - C/2 \\
&\geq& \sum_{ij}C_{ij}(I^*(t), S(t))\hat{W}_{ij}(t) - C/2 \\
&\geq& \sum_{ij}C_{ij}(I^*(t), S(t)) W_{ij}(t) - C
\end{eqnarray*}
It follows that the resource allocation $\hat{I}(t)$ (and the corresponding 
transmission decisions given by (\ref{eq:mu-ij-enhanced})) yields a $C$-approximation.
\end{proof}

\section{Approximate Scheduling and Slater Conditions} \label{section:slater}

Here we replace Assumption A2 with a stronger assumption that states
the constraints can be satisfied with $\delta$ slackness. This is related to a \emph{Slater condition} 
in classical static optimization problems \cite{bertsekas-nonlinear}.  It allows all queues to be
deterministically bounded.  It also allows performance analysis for implementations when the 
error in the 
attempted 
minimization of the right hand side of (\ref{eq:big-1-drift}) is off by more than just a constant $C$, 
such as an amount that
may be proportional to the queue backlog  (similar to the $\theta$-multiplicative-approximations discussed
in Section \ref{section:internet-wireless}). 

\emph{Assumption A3}: There exists a value $\delta>0$ such that for all 
$\omega \in \{\omega(0), \ldots, \omega(t_{end} - 1)\}$, there is
at least one control action $\alpha_{\omega}' \in \script{A}_{\omega}$ that satisfies: 
\begin{eqnarray*}
&\hat{y}_l(\alpha_{\omega}', \omega) + g_l(\hat{\bv{x}}(\alpha'_{\omega}, \omega))  \leq -\delta \: \: \forall l \in \{1, \ldots, L\}& \\
&\hat{a}_k(\alpha_{\omega}', \omega) \leq \hat{b}_k(\alpha_{\omega}', \omega) - \delta \: \: \forall k \in \{1, \ldots, K\}& \\
&\hat{\bv{x}}(\alpha_{\omega}', \omega) + \bv{\epsilon} \in \script{X}& 
\end{eqnarray*}
for all vectors $\bv{\epsilon} = (\epsilon_1, \ldots, \epsilon_M)$ with entries $\epsilon_m$
that satisfy $|\epsilon_m| \leq \delta$ for all $m \in\{1, \ldots, M\}$.  Further, assume that: 
\[ x_m^{min} + \delta \leq \hat{x}_m(\alpha_{\omega}', \omega) \leq x_m^{max} - \delta \: \: \forall m \in \{1, \ldots, M\} \]
This final assumption is mild and can easily be engineered to be true
by convexly extending the range of the convex functions
$f(\cdot)$, $g_l(\cdot)$ by $\delta$ in all directions (so that $x_m^{min}$ is decreased by $\delta$ and $x_m^{max}$
is increased by $\delta$), as in \cite{neely-delay-based-infocom2010}. 

Define a \emph{$C(t)$-approximation} as an algorithm that, every slot $t$, observes the current 
queue states and 
choses a control action that comes within $C(t)$ of minimizing the right 
hand side of (\ref{eq:big-1-drift}), where $C(t)$ is a value that can depend on $t$.
Suppose that we implement the universal scheduling algorithm of Section \ref{section:general-alg} 
using a $C(t)$-approximation with $C(t)$ that satisfies the following for all $t$: 
\begin{eqnarray}
C(t) &\leq&  C + V\epsilon_V   +  \sum_{l=1}^L Z_l(t)\epsilon_Z \nonumber  \\
&& + \sum_{k=1}^KQ_k(T)\epsilon_Q 
+ \sum_{m=1}^M|H_m(t)|\epsilon_H \label{eq:ct}
\end{eqnarray}
where $C$, $\epsilon_V$, $\epsilon_Z$, $\epsilon_Q$, $\epsilon_H$ are non-negative constants. 
Note that this is a $C$-approximation if $\epsilon_V = \epsilon_Z = \epsilon_Q = \epsilon_H = 0$, and is 
the exact minimization of (\ref{eq:big-1-drift}) if we additionally have $C=0$. 

\begin{thm}  \label{thm:slater} Suppose 
Assumptions A1 and A3 hold for
some $\delta>0$. Consider any $C(t)$-approximate algorithm that satisfies 
(\ref{eq:ct}) every slot $t$, and assume that: 
\begin{eqnarray*}
 \epsilon_Q < \delta \: \: , \: \:  \epsilon_H < \delta \\
\epsilon_Z + \epsilon_H \sum_{m=1}^M \beta_{l,m} < \delta \: \: \forall l \in \{1, \ldots, L\}
\end{eqnarray*}
Let the random event sequence 
$\{\omega(0), \omega(1), \ldots \}$ be arbitrary.  Suppose all initial queue backlogs are zero. 
Then: 

(a) All queue backlogs are bounded, so that for any slot $t\geq 0$ we have: 
\[ Q_k(t), Z_l(t), |H_m(t)| \leq \frac{VC_3}{\theta} \]
where 
\[ C_3 \defequiv \sqrt{D_1 + D_2 + D_3} \]
where $D_1$, $D_2$, $D_3$ are constants defined as: 
\begin{eqnarray*} 
D_1 &\defequiv& \left[\frac{B+C}{V} + (y_0^{max} - y_0^{min}) + (f^{max} - f^{min}) + \epsilon_V\right]^2  \\
D_2 &\defequiv&  2D\theta^2/V^2 \\
D_3 &\defequiv& \frac{2z_{max}\theta}{V}\sqrt{D_1} 
\end{eqnarray*}
where $D$ is defined in (\ref{eq:D}), 
$z_{max}$ is the maximum over all $z_l^{diff}$, 
$q_k^{diff}$, and $h_m^{diff}$ constants, and 
where $\theta$ is defined: 
\begin{equation} \label{eq:theta} 
 \theta \defequiv \min\left[\delta - \epsilon_Q, \delta - \epsilon_H, \frac{\delta - \epsilon_Z - \epsilon_H\beta_{sum}}{1 + \beta_{sum}}\right] 
 \end{equation} 
 where $\beta_{sum}$ is defined: 
 \[ \beta_{sum}\defequiv \max_{l \in \{1, \ldots, L\}} \sum_{m=1}^M\beta_{l,m} \]
  In the special case when $\epsilon_Z = \epsilon_Q = \epsilon_H = 0$, we have
  $\theta = \delta/(1+\beta_{sum})$.

(b) For any designated time $t_{end}>0$ we have: 
\begin{eqnarray*}
\overline{y}_l + g_l(\overline{\bv{x}}) \leq  \frac{VC_3}{\theta t_{end}}\left[1 + \sum_{m=1}^M\beta_{l,m}\right] &  \forall l \in\{1, \ldots, L\}\\
\overline{a}_k \leq \overline{b}_k + \frac{VC_3}{\theta t_{end}} &  \forall k \in \{1, \ldots, K\} \\
\overline{\bv{x}} + \overline{\bv{\epsilon}}(t) \in \script{X} & 
\end{eqnarray*}
where $\bv{\epsilon}(t) = (\epsilon_1(t), \ldots, \epsilon_M(t))$ has entries that satisfy: 
\[ |\epsilon_m(t)| \leq  \frac{VC_3}{\theta t_{end}} \]

c) Consider any positive integer frame size $T$, any positive integer $R$, and define
$t_{end} = RT$.  Then the value of the system cost metric over $t_{end}$ slots satisfies: 
\begin{eqnarray}
\overline{y}_0 + f(\overline{\bv{x}}) \leq \frac{(1-p)}{R}\sum_{r=0}^{R-1} F_r^* + \epsilon_V + p(y_0^{max} + f^{max}) \nonumber \\
+ \frac{B + C + \tilde{D}(T-1)}{V} +  \sum_{m=1}^M\frac{\nu_mVC_3}{\theta RT}  \label{eq:util-slater} 
\end{eqnarray}
where $p$ is defined: 
\[ p \defequiv \max\left[\frac{\epsilon_Z}{\delta - (\epsilon_H+\theta)\beta_{sum}}, \frac{\epsilon_Q}{\delta}, 
\frac{\epsilon_H}{\epsilon_H + \theta}\right] \]
and where $\tilde{D}$ is a constant defined in (\ref{eq:d-tilde}), $B$ is a constant 
defined in (\ref{eq:B}), and $C_3$ is a constant defined
in part (a). 
\end{thm} 
\begin{proof} 
See Appendix D. 
\end{proof} 

The cost bound (\ref{eq:util-slater}) can be understood as follows: The last term on the right hand side
goes to zero as $R$ increases (and is equal to $0$ for all $R$ if $f(\cdot) = 0$ so that $\nu_m = 0$ for all $m$). 
The second to last term can be made arbitrarily small with a suitably large $V$. 
Finally, if $\epsilon_Z$, $\epsilon_Q$, $\epsilon_H$, $\epsilon_V$ are small, then $p$ is small and the remaining
terms on the right hand side are close to the cost $\frac{1}{R}\sum_{r=0}^{R-1}F_r^*$ which is associated with 
implementing the $T$-slot lookahead policy over $R$ frames.

\section{Conclusions} 

We have developed a framework for universal constrained optimization of time averages in 
time varying systems.  Our results hold for any event sample paths and do not require a probability 
model. It was shown that performance can closely track the performance of an ideal policy with 
knowledge of the future up to $T$ slots, provided that we allow the number of $T$-slot frames, denoted by $R$, 
to be large enough to ensure that the error
terms decay to a negligible value.  This framework 
was applied to an internet model and to a more extensive
queueing network model to provide utility guarantees with deterministic queue bounds for
arbitrary traffic, channels, and mobility. 


\section*{Appendix A --- Proof of Theorem \ref{thm:1}} 

We first prove parts (a) and (b) of Theorem \ref{thm:1}. 

\begin{proof} (Theorem \ref{thm:1} part (a)) Let $\bv{\gamma}(t)$ and $\alpha(t)$ represent the decisions
made by the $C$-approximate policy on slot $t$, which necessarily 
satisfy constraints (\ref{eq:new-c4})-(\ref{eq:new-c6}).  
Because these decisions come within $C$ of minimizing
the right hand side of (\ref{eq:big-1-drift}) over all other possible decisions, 
we have from (\ref{eq:big-1-drift}): 
\begin{eqnarray} 
\Delta_1(t) + V\hat{y}_0(\alpha(t), \omega(t)) + Vf(\bv{\gamma}(t)) \leq \nonumber \\
B + C + V\hat{y}_0(\alpha^*(t), \omega(t)) + Vf(\bv{\gamma}^*(t)) \nonumber \\
+ \sum_{l=1}^LZ_l(t)[\hat{y}_l(\alpha^*(t), \omega(t)) + g_l(\bv{\gamma}^*(t))] \nonumber \\
+  \sum_{k=1}^KQ_k(t)[\hat{a}_k(\alpha^*(t), \omega(t)) - \hat{b}_k(\alpha^*(t), \omega(t))]  \nonumber \\
+ \sum_{m=1}^MH_m(t)[\gamma_m^*(t) - \hat{x}_m(\alpha^*(t), \omega(t))]  \label{eq:big2} 
\end{eqnarray} 
where $\alpha^*(t)$ and $\bv{\gamma}^*(t)$ represent any alternative decisions that could be made
on slot $t$ that satisfy (\ref{eq:new-c4})-(\ref{eq:new-c6}). 

Now choose $\alpha^*(t) = \alpha'_{\omega(t)}$, where $\alpha'_{\omega(t)} \in \script{A}_{\omega(t)}$ 
is the decision known to exist by 
Assumption A2 that satisfies:
\begin{eqnarray*} 
& \hat{y}_l(\alpha^*(t), \omega(t)) + g_l(\hat{\bv{x}}(\alpha^*(t), \omega(t))) \leq 0 &  \forall l \in \{1, \ldots, L\} \\
&  \hat{a}_k(\alpha^*(t), \omega(t)) \leq \hat{b}_k(\alpha^*(t), \omega(t)) &  \forall k \in \{1, \ldots, K\}  \\
&  \hat{\bv{x}}(\alpha^*(t), \omega(t)) \in \script{X} & 
\end{eqnarray*}
Further, choose 
$\bv{\gamma}^*(t) = \hat{\bv{x}}(\alpha^*(t), \omega(t))$. These decisions satisfy (\ref{eq:new-c4})-(\ref{eq:new-c6}), and  
plugging these decisions directly into the right hand side of (\ref{eq:big2}) yields: 
\begin{eqnarray*} 
\Delta_1(t) + V\hat{y}_0(\alpha(t), \omega(t)) + Vf(\bv{\gamma}(t)) \leq \nonumber \\
B + C  + V\hat{y}_0(\alpha^*(t), \omega(t)) + Vf(\bv{\gamma}^*(t)) 
\end{eqnarray*} 
Rearranging terms and using the bounds $y_0^{min}$, $y_0^{max}$ and 
$f^{min}$, $f^{max}$ yields: 
\[ \Delta_1(t) \leq B + C + V(y_0^{max} - y_0^{min}) + V(f^{max} - f^{min}) \]
Let the right hand side of the above inequality be denoted by $P$. Using the 
definition of $\Delta_1(t)$ thus gives: 
\[ L(\bv{\Theta}(t+1)) - L(\bv{\Theta}(t)) \leq P  \]
The above holds for all $t\geq0$.  Summing
over $\tau \in \{0, \ldots, t-1\}$ (for some time $t>0$) and dividing by $t$ yields: 
\[ \frac{1}{t}[L(\bv{\Theta}(t)) - L(\bv{\Theta}(0))] \leq  P  \]
Using the definition of $L(\bv{\Theta}(t))$ in (\ref{eq:L})  
proves part (a). 
\end{proof} 

\begin{proof} (Theorem \ref{thm:1} part (b)) From part (a), 
if all queues are initially empty (so that $L(\bv{\Theta}(0)) = 0$), we have for all 
slots $t>0$: 
\begin{equation} \label{eq:plug-ineq} 
Z_l(t) , Q_k(t), |H_m(t)|  \leq C_0\sqrt{tV}  
\end{equation} 
Plugging (\ref{eq:plug-ineq}) into (\ref{eq:foo2}), (\ref{eq:truesat-g}), (\ref{eq:truesat-x})  of Lemma \ref{lem:constraints}
proves part (b).
\end{proof} 

To prove part (c) of Theorem \ref{thm:1}, we need the following preliminary lemma. 

\begin{lem} \label{lem:T-slot-approx} For any initial time $t_0$, any queue values $\bv{\Theta}(t_0)$, 
any integer $T>1$, and any collection of
$C$-approximate decisions that are implemented over the $T$-slot interval $\tau \in \{t_0, \ldots, t_0 + T-1\}$, 
we have: 
\begin{eqnarray*}
\Delta_T(t_0) + \sum_{\tau=t_0}^{t_0+T-1} [V\hat{y}_0(\alpha(\tau), \omega(\tau)) + Vf(\bv{\gamma}(\tau)) ] \leq \\
BT + CT +  DT(T-1) \\
+ \sum_{\tau=t_0}^{t_0+T-1} [V\hat{y}_0(\alpha^*(\tau), \omega(\tau)) + Vf(\bv{\gamma}^*(\tau)) ]  \\
+ \sum_{l=1}^LZ_l(t_0) \sum_{\tau=t_0}^{t_0+T-1}\left[  \hat{y}_l(\alpha^*(\tau), \omega(\tau)) +  g_l(\bv{\gamma}^*(\tau)) \right] \\
+ \sum_{k=1}^KQ_k(t_0)\sum_{\tau=t_0}^{t_0+T-1} \left[\hat{a}_k(\alpha^*(\tau), \omega(\tau)) - \hat{b}_k(\alpha^*(\tau), \omega(\tau))\right]  \\
+ \sum_{m=1}^MH_m(t_0) \sum_{\tau=t_0}^{t_0+T-1} \left[\gamma_m^*(\tau) - \hat{x}_m(\alpha^*(\tau), \omega(\tau))\right]
\end{eqnarray*}
for any alternative decisions $\alpha^*(\tau)$, $\bv{\gamma}^*(\tau)$ over $\tau \in \{t_0, \ldots, t_0 + T-1\}$ that
satisfy (\ref{eq:new-c4})-(\ref{eq:new-c6}).  The constant $B$ is defined according to (\ref{eq:B}) and 
$D$ is defined by (\ref{eq:D}).
\end{lem} 

\begin{proof} (Lemma \ref{lem:T-slot-approx})
Because our policy is $C$-approximate, for all slots $t$ it comes within $C$ of minimizing
the right hand side of (\ref{eq:big-1-drift}).  Hence, for all $\tau \in \{t_0, \ldots, t_0+T-1\}$ we have: 
\begin{eqnarray} 
\Delta_1(\tau) + V\hat{y}_0(\alpha(\tau), \omega(\tau)) + Vf(\bv{\gamma}(\tau))  \leq \nonumber \\
B + C + V\hat{y}_0(\alpha^*(\tau), \omega(\tau)) + Vf(\bv{\gamma}^*(\tau)) \nonumber \\
+ \sum_{l=1}^LZ_l(\tau)[\hat{y}_l(\alpha^*(\tau), \omega(\tau))  + g_l(\bv{\gamma}^*(\tau))] \nonumber \\
+  \sum_{k=1}^KQ_k(\tau)[\hat{a}_k(\alpha^*(\tau), \omega(\tau)) - \hat{b}_k(\alpha^*(\tau), \omega(\tau))] \nonumber \\
+ \sum_{m=1}^MH_m(\tau)[\gamma_m^*(\tau) - \hat{x}_m(\alpha^*(\tau), \omega(\tau))]\label{eq:T-slot-approx1} 
\end{eqnarray} 
However, by definition of $z_l^{diff}$, $q_k^{diff}$, $h_m^{diff}$, the queues
$Z_l(t)$, $Q_k(t)$, $H_m(t)$ can change by at most these values
on each slot, and hence for $\tau \in \{t_0, \ldots, t_0 + T-1\}$ we have: 
\begin{eqnarray*} 
|Z_l(\tau) - Z_l(t_0)| &\leq& z_l^{diff}\cdot(\tau - t_0) \\
|Q_k(\tau) - Q_k(t_0)| &\leq& q_k^{diff}\cdot(\tau-t_0) \\
|H_m(\tau) - H_m(t_0)| &\leq& h_m^{diff}\cdot(\tau-t_0)
\end{eqnarray*}
Using these in (\ref{eq:T-slot-approx1}) gives: 
\begin{eqnarray*} 
\Delta_1(\tau) + V\hat{y}_0(\alpha(\tau), \omega(\tau)) + Vf(\bv{\gamma}(\tau))  \leq \nonumber \\
B + C + 2D\cdot(\tau - t_0) + V\hat{y}_0(\alpha^*(\tau), \omega(\tau)) + Vf(\bv{\gamma}^*(\tau)) \nonumber \\
+ \sum_{l=1}^LZ_l(t_0)[\hat{y}_l(\alpha^*(\tau), \omega(\tau))  + g_l(\bv{\gamma}^*(\tau))] \nonumber \\
+  \sum_{k=1}^KQ_k(t_0)[\hat{a}_k(\alpha^*(\tau), \omega(\tau)) - \hat{b}_k(\alpha^*(\tau), \omega(\tau))] \nonumber \\
+ \sum_{m=1}^MH_m(t_0)[\gamma_m^*(\tau) - \hat{x}_m(\alpha^*(\tau), \omega(\tau))]
\end{eqnarray*} 
where $D$ is defined in (\ref{eq:D}). Summing the above inequality over $\tau \in \{t_0, \ldots, t_0+T-1\}$
yields the result, where we use the fact that: 
\[ \sum_{\tau=t_0}^{t_0+T-1} (\tau-t_0) = T(T-1)/2 \]
\end{proof}

We can now prove Theorem \ref{thm:1} part (c). 

\begin{proof} (Theorem \ref{thm:1} part (c)) Fix integers $r \geq 0$ and $T>0$. 
Fix $\epsilon>0$, and 
let $\alpha^*(\tau)$ represent the
decisions over the interval $\tau \in \{rT, \ldots, (r+1)T - 1\}$ that solve the 
problem (\ref{eq:frame-lookahead-problem}) and yield cost that is no 
more than $F_r^* + \epsilon$. 
Let $\bv{\gamma}^*(\tau)$ be constant over $\tau \in \{rT, \ldots, (r+1)T-1\}$, given by: 
\[ \bv{\gamma}^*(\tau) = \frac{1}{T}\sum_{t=rT}^{(r+1)T-1} \hat{\bv{x}}(\alpha^*(t), \omega(t)) \]
Plugging these alternative decisions $\alpha^*(\tau)$ and $\bv{\gamma}^*(\tau)$ 
into the result of Lemma \ref{lem:T-slot-approx} for $t_0 = rT$
yields: 
\begin{eqnarray*}
\Delta_T(rT) + \sum_{\tau=rT}^{(r+1)T-1} [V\hat{y}_0(\alpha(\tau), \omega(\tau)) + Vf(\bv{\gamma}(\tau)) ] \leq \\
BT +  CT +  DT(T-1) 
+VT(F^*_r + \epsilon)  
\end{eqnarray*}
The above holds for all $\epsilon>0$, and hence we can take a limit as $\epsilon\rightarrow 0$ to remove the 
$\epsilon$ in the final 
term. 
Define $t_{end} = RT$ for some positive integer $R$.  Summing the above over $r \in \{0, \ldots, R-1\}$
and dividing by $VRT$ yields: 
\begin{eqnarray}
\overline{y}_0 + f(\overline{\bv{\gamma}}) + \frac{L(\bv{\Theta}(RT)) - L(\bv{\Theta}(0))}{VRT} \leq \nonumber \\
\frac{1}{R}\sum_{r=0}^{R-1}F_r^* + \frac{B + C}{V}  +  \frac{D(T-1)}{V} \label{eq:last-part-proof-1c}
\end{eqnarray}
where $\overline{y}_0$ and $\overline{\bv{\gamma}}$
represent time averages over the first $t_{end}$ slots, 
and where we have used Jensen's inequality in the concave function $f(\bv{\gamma})$.

However, we have by (\ref{eq:cost-assumption-f}): 
\begin{eqnarray*} 
f(\overline{\bv{x}}) &\leq& f(\overline{\bv{\gamma}}) + \sum_{m=1}^M\nu_m |\overline{\gamma}_m - \overline{x}_m| \\
&=& f(\overline{\bv{\gamma}}) + \sum_{m=1}^M\frac{\nu_m|H_m(RT)-H_m(0)|}{RT}
\end{eqnarray*}
where the final equality holds by (\ref{eq:foo3}). 
Using this in (\ref{eq:last-part-proof-1c})  together with the fact that
$L(\cdot) \geq 0$ yields the cost bound (\ref{eq:utility-thm}) of part (c).  
Finally, if initial queue backlogs are 0, by (\ref{eq:plug-ineq}) applied to time $t=RT$ we have: 
\[ \frac{|H_m(RT) - H_m(0)|}{RT} = \frac{|H_m(RT)|}{RT} \leq \frac{C_0\sqrt{V}}{\sqrt{RT}} \]
\end{proof}

\section*{Appendix B --- Conditions for Achievability of $F^*$} 

Consider the following additional assumption. 

\emph{Assumption A4:}  We have either one of the following two conditions: 
\begin{enumerate} 
\item For all $\omega \in \{\omega(0), \ldots, \omega(t_{end}-1)\}$, the control action 
space $\script{A}_{\omega}$ contains a finite number of actions. 

\item For all $\omega \in \{\omega(0), \ldots, \omega(t_{end}-1)\}$, 
the set $\script{A}_{\omega}$ is a compact subset of $\mathbb{R}^c$ for some dimension $c$, 
the functions $\hat{y}_l(\alpha, \omega), \hat{a}_k(\alpha, \omega)$ are lower semi-continuous over $\alpha \in  \script{A}_{\omega}$, the
functions $\hat{b}_k(\alpha, \omega)$ are upper semi-continuous over $\alpha \in  \script{A}_{\omega}$, and 
the functions $\hat{x}_m(\alpha, \omega)$ are continuous over $\alpha \in \script{A}_{\omega}$.   Note that
all continuous functions are both upper and lower semi-continuous.\footnote{A function $b(\bv{\alpha})$ is \emph{upper semi-continuous} over $\bv{\alpha} \in \script{A}$ if for any $\bv{\alpha} \in \script{A}$, we 
have $b(\bv{\alpha}) \geq \lim_{n\rightarrow\infty} b(\bv{\beta}_n)$ for all sequences $\bv{\beta}_n \in \script{A}$ such that 
$\lim_{n\rightarrow\infty} \bv{\beta}_n = \bv{\alpha}$.  A function is \emph{lower semi-continuous} if the inequality is 
reversed.  All  bounded functions that are discontinuous only on a set of measure zero can be 
easily modified to have the desired semi-continuous property by appropriately
re-defining the function value at  points of discontinuity. Most systems of practical interest 
have the desired semi-continuity
properties.} 
\end{enumerate}

\begin{lem} \label{lem:existence}  Suppose Assumptions A1, A2, A4 hold for 
given values $\{\omega(0), \ldots, \omega(t_{end}-1)\}$.  Then the infimum value $F^*$ for the 
problem (\ref{eq:opt})-(\ref{eq:c4}) can be 
achieved by a particular (possibly non-unique) 
sequence of control actions $\{\alpha^*(0), \ldots, \alpha^*(t_{end}-1)\}$.  
That is, these actions satisfy the feasibility constraints (\ref{eq:c1})-(\ref{eq:c4}), and yield: 
\[ \overline{y}_0 + f(\overline{x}_1, \ldots, \overline{x}_M) = F^* \]
where: 
\begin{eqnarray*}
\overline{x}_m &=& \frac{1}{t_{end}} \sum_{\tau=0}^{t_{end}-1} \hat{x}_m(\alpha^*(\tau), \omega(\tau)) \: \: \: \:  \forall m \in \{1, \ldots, M\} 
\end{eqnarray*}
and where $\overline{y}_0$ is similarly defined as a time average over $\tau \in \{0, \ldots, t_{end}-1\}$. 
\end{lem} 

\begin{proof} 
We already know that Assumption A2 implies the existence of a feasible sequence of control actions. 
Thus, the infimum value $F^*$ of the cost metric over all feasible policies 
is well defined, and by Assumption A1 it must satisfy: 
\[ y_0^{min} + f^{min} \leq F^* \leq y_0^{max} + f^{max} \]
Consider now the case when  the first condition of Assumption A4 holds.  Then there are only a finite
number of possible control sequences over the horizon $\{0, 1, \ldots, t_{end} - 1\}$, and
so there is one that achieves the minimum cost  value $F^*$.  The case when the second
condition of Assumption A4 holds can be proven using the Bolzano-Wierstrass Theorem together
with a simple limiting argument, and is omitted for brevity.
\end{proof}

\section*{Appendix C -- Ergodicity} 

Consider the infinite horizon problem discussed in Section \ref{section:infinite-horizon}. 
Suppose that the random events $\{\omega(0), \omega(1), \omega(2), \ldots\}$ evolve according to a 
general ergodic process with a well defined time average probability distribution.  Specifically, let $\Omega$ represent 
a finite (but arbitrarily large) outcome space for $\omega(t)$, and for each $\omega \in \Omega$ assume that
there is a steady state value $\pi(\omega)$, such that: 
\begin{eqnarray*}
\lim_{t\rightarrow\infty} \frac{1}{t}\sum_{\tau=0}^{t-1} 1_{\omega}(\tau) = \pi(\omega) \: \: \mbox{ with probability 1} 
\end{eqnarray*}
where $1_{\omega}(\tau)$ is an indicator function that is $1$ if $\omega(\tau) = \omega$, and zero else.  Further,
assume the limiting probability converges uniformly to the steady state value, regardless of past history, so that: 
\begin{eqnarray*}
|Pr[\omega(t+ t_0) = \omega | History(t_0)] - \pi(\omega)| \leq error(t)
\end{eqnarray*}
where $History(t_0)$ represents the past history of the process up to slot $t_0$, and 
where $error(t)$ is a function that decays to $0$ as $t\rightarrow \infty$, regardless of the past history.
This is related to the \emph{decaying memory property} in \cite{neely-stock-arxiv} and the admissibility
assumptions in \cite{neely-thesis}\cite{now}.  

In this case, it can be shown that the 
optimal infinite horizon cost, denoted $f^*$, can be achieved over the class of stationary
and randomized algorithms that make (possibly probabilistic) decisions for control actions on each slot $t$
based only on the 
current state $\omega(t)$ (see \cite{neely-thesis}\cite{neely-power-network-jsac}\cite{neely-energy-it} for related
proofs).\footnote{Similar results on optimality of stationary policies can typically be 
achieved when the cardinality of the set $\Omega$ 
is infinite, although steady state time averages and uniform convergence to a steady state
are more awkward to deal with in this case.  The easiest such arguments for (possibly uncountably) infinite
sets $\Omega$ are for $\omega(t)$ processes that are i.i.d. over slots.}  Further, under mild conditions (such as the existence
of a value $\delta>0$ for which the Slater condition of Assumption A3 in Section \ref{section:slater} 
is satisfied), the value of $F_r^*$, being the 
optimal time average 
cost under the $T$-slot lookahead policy over the $T$-slot interval starting at time $rT$, satisfies for any integer $r \geq0$: 
\[ \lim_{T\rightarrow\infty} F_r^* = f^* \: \: \mbox{ with probability 1} \]
That is, regardless of the past history before time $rT$, the $T$-slot lookahead policy over a very large $T$ approaches 
the optimal $f^*$.   The reason the ``mild'' additional conditions, such as the Slater condition, is needed, is that
$F_r^*$ requires all constraints to be exactly satisfied by the end of the $T$ slots,  whereas the infinite horizon problem
does not require this.  

Because of the uniform error decay, we have: 
\[ |\expect{F_r^*} -  f^*| \leq   \delta(T) \]
where $\delta(T)$ is a function such that $\delta(T) \rightarrow 0$ as $T\rightarrow \infty$. 
Therefore we can write: 
\[ F_r^* = f^* + \delta_r \]
where $\delta_r$ is a random variable that satisfies $|\expect{\delta_r}| \leq \delta(T)$ for all $r$. 

Using the definition of $\delta_r$, it follows from (\ref{eq:infinite-horizon}) that time average cost
satisfies for any integer $T>0$: 
\begin{eqnarray*}
 \lim_{R\rightarrow\infty} [\overline{y}_0(RT) + f(\overline{\bv{x}}(RT))] &\leq& f^* + \frac{B+C}{V}  + \frac{D(T-1)}{V} \\
&& +  \lim_{R\rightarrow\infty}\frac{1}{R}\sum_{r=0}^{R-1} \delta_r  
 \end{eqnarray*}
Under mild conditions, such as when $\{\omega(t)\}$ evolves according to a finite state ergodic
 Markov chain, 
law of large number averaging principles imply that the last term, being a time average of the $\delta_r$
values, is bounded in absolute value with probability 1 by $\delta(T)$, a term that is negligibly small for
large values of $T$.  We can choose a large $T$ provided that we also compensate with a large $V$
to make the $D(T-1)/V$ term neglibible. 
This demonstrates that, with probability 1, the algorithm implemented over an 
infinite time horizon yields cost that can be pushed arbitrarily close to the optimal value $f^*$ if
$V$ is suitably large. 

\section*{Appendix D -- Proof of Theorem \ref{thm:slater}} 

\begin{proof} (Theorem \ref{thm:slater}a)  
Define $\theta$ as the positive real number that solves the following problem (it can be shown that
the solution is given by (\ref{eq:theta})): 
\begin{eqnarray}
\mbox{Maximize:} & \theta   \label{eq:theta-problem} \\
\mbox{Subject to:} &\theta \leq \delta - \epsilon_Q \nonumber \\
&\theta \leq \delta - \epsilon_Z - (\epsilon_H + \theta)\sum_{m=1}^M\beta_{l,m} \: \: \forall l \nonumber  \\
&  \theta \leq \delta - \epsilon_H \nonumber
\end{eqnarray}
Following the proof of Theorem \ref{thm:1}a and replacing $C$ with $C(t)$ 
we have (compare with (\ref{eq:big2})): 
\begin{eqnarray} 
\Delta_1(t) + V\hat{y}_0(\alpha(t), \omega(t)) + Vf(\bv{\gamma}(t)) \leq \nonumber \\
B + C + V\hat{y}_0(\alpha^*(t), \omega(t)) + Vf(\bv{\gamma}^*(t)) + V\epsilon_V\nonumber \\
+ \sum_{l=1}^LZ_l(t)[\epsilon_Z + \hat{y}_l(\alpha^*(t), \omega(t)) + g_l(\bv{\gamma}^*(t))] \nonumber \\
+  \sum_{k=1}^KQ_k(t)[\epsilon_Q  + \hat{b}_k(\alpha^*(t), \omega(t)) - \hat{a}_k(\alpha^*(t), \omega(t))]  \nonumber \\
+ \sum_{m=1}^M|H_m(t)|\epsilon_H  \nonumber \\
+ \sum_{m=1}^MH_m(t)[\gamma_m^*(t) - \hat{x}_m(\alpha^*(t), \omega(t))]  \label{eq:biga3} 
\end{eqnarray} 
where $\alpha^*(t)$ and $\bv{\gamma}^*(t)$ represent any alternative decisions that could be made
on slot $t$ that satisfy (\ref{eq:new-c4})-(\ref{eq:new-c6}). 

Now choose $\alpha^*(t) = \alpha'_{\omega(t)}$, where $\alpha'_{\omega(t)} \in \script{A}_{\omega(t)}$ is the 
decision known to exist by Assumption A3.  Further, choose $\bv{\gamma}^*(t) = \bv{\gamma}'(t)$, where
$\bv{\gamma}'(t) = (\gamma_m'(t))|_{m=1}^M$ is defined such that for all $m \in \{1, \ldots, M\}$: 
\begin{equation} \label{eq:gamma-prime} \gamma_m'(t) \defequiv \left\{ \begin{array}{ll}
                          \hat{x}_m(\alpha'_{\omega(t)}, \omega(t)) - \epsilon_H - \theta &\mbox{ if $H_m(t)\geq0$} \\
                              \hat{x}_m(\alpha'_{\omega(t)}, \omega(t)) + \epsilon_H + \theta  & \mbox{ if $H_m(t) < 0$} 
                            \end{array}
                                 \right.
\end{equation} 
This is feasible because of the last inequality in Assumption A3 together with the fact that: 
\begin{eqnarray*}
 -\delta \leq -\epsilon_H - \theta \leq \epsilon_H  + \theta \leq \delta
 \end{eqnarray*}
With these choices, (\ref{eq:biga3}) becomes: 
\begin{eqnarray} 
\Delta_1(t)  \leq \nonumber \\
B + C + V(y_0^{max} - y_0^{min})+ V(f^{max} - f^{min})+ V\epsilon_V\nonumber \\
- \sum_{l=1}^LZ_l(t)[\delta - \epsilon_Z - (\epsilon_H+\theta)\sum_{m=1}^M\beta_{l,m}] \nonumber \\
-  \sum_{k=1}^KQ_k(t)[\delta - \epsilon_Q]  \nonumber \\
-\sum_{m=1}^M|H_m(t)|\theta \label{eq:biga4}  
\end{eqnarray} 
where have used the fact that, from (\ref{eq:cost-assumption-g}): 
\begin{eqnarray*}
g_l(\bv{\gamma}^*(t)) \leq g_l(\hat{\bv{x}}(\alpha'_{\omega(t)}, \omega(t))) + (\epsilon_H + \theta)\sum_{m=1}^M\beta_{l,m}
\end{eqnarray*}
Now define $P$ as: 
\[ P \defequiv B + C + V(y_0^{max} - y_0^{min})+ V(f^{max} - f^{min})+ V\epsilon_V \]
Because 
the value $\theta$ is a bound on all the terms multiplying queue values in (\ref{eq:biga4}), we have: 
\begin{eqnarray} 
\Delta_1(t)  \leq P
- \theta\sum_{l=1}^LZ_l(t)
-  \theta\sum_{k=1}^KQ_k(t) - \theta \sum_{m=1}^M|H_m(t)| \label{eq:biga5}  
\end{eqnarray} 
It follows that the drift is non-positive whenever the sum of the absolute value of queue size is greater than or equal 
to $P/\theta$. It is not difficult to show that the largest possible value of $L(\bv{\Theta}(t))$ under the constraint
that the sum of absolute queue values is less than or equal to $P/\theta$ is $(1/2)(P/\theta)^2$.  Hence, if the Lyapunov function 
is larger than this value, it cannot increase on the next slot.  However, if the absolute sum on slot $t$ is less than 
or equal to $P/\theta$ we have: 
\begin{eqnarray}
L(\bv{\Theta}(t+1)) &\leq& \frac{1}{2}\sum_{l=1}^L(Z_l(t) + z_l^{diff}(t))^2 \nonumber \\
&& + \frac{1}{2}\sum_{k=1}^K(Q_k(t) + q_k^{diff}(t))^2 \nonumber \\
&& + \frac{1}{2}\sum_{m=1}^M(|H_m(t)|  + h_m^{diff}(t))^2 \nonumber \\
&\leq& L(\bv{\Theta}(t)) + D + z_{max}P/\theta \nonumber \\
&\leq& (1/2)(P/\theta)^2 + D + z_{max}P/\theta \label{eq:calc1} 
\end{eqnarray}
where $z_l^{diff}(t)$, $q_k^{diff}(t)$, $h_m^{diff}(t)$
represent the absolute value of the 
change in $Z_l(t)$, $Q_k(t)$, $H_m(t)$, respectively, over one slot, having
maximum absolute value given by $z_l^{diff}$, 
$q_k^{diff}$, and $h_m^{diff}$,  $D$ is defined in (\ref{eq:D}), and
$z_{max}$ is the maximum over all $z_l^{diff}$, 
$q_k^{diff}$, and $h_m^{diff}$ constants.

It follows that for all $t$ we have:\footnote{More precisely, the bound on $L(\bv{\Theta}(t))$ 
is clearly true for $t=0$.  Supposing
it is true for slot $t$, we show it is true for slot $t+1$:  If the absolute sum is greater than or equal to  
$P/\theta$ on slot $t$, then the Lyapunov value cannot increase on the next slot and so the bound
also holds for slot $t+1$.  Else, if the absolute sum is less than $P/\theta$ on slot $t$, 
then the bound again holds for slot $t+1$ by the calculation (\ref{eq:calc1}).}
\[ L(\bv{\Theta}(t)) \leq  (1/2)(P/\theta)^2 + D + z_{max}P/\theta \]
Therefore, all queues are bounded by: 
\[ Q_k(t), Z_l(t), |H_m(t)| \leq \sqrt{(P/\theta)^2 + 2D + 2z_{max}P/\theta }  \]
This bound is given by: 
\begin{eqnarray*}
\frac{V\sqrt{ P^2/V^2+ 2D\theta^2/V^2 + 2z_{max}\theta P/V^2}}{\theta} \\
= \frac{V\sqrt{D_1 + D_2 + D_3}}{\theta} 
\end{eqnarray*} 
where 
\begin{eqnarray*} 
D_1 &\defequiv& \left[\frac{B+C}{V} + (y_0^{max} - y_0^{min}) + (f^{max} - f^{min}) + \epsilon_V\right]^2  \\
D_2 &\defequiv&  2D\theta^2/V^2 \\
D_3 &\defequiv& \frac{2z_{max}\theta}{V}\sqrt{D_1} 
\end{eqnarray*}
\end{proof} 

\begin{proof} (Theorem \ref{thm:slater}b)
The proof follows immediately 
by applying the queue bounds of part (a) to the constraint 
bounds (\ref{eq:foo2}), (\ref{eq:truesat-g}), (\ref{eq:truesat-x}) 
of Lemma \ref{lem:constraints}, using initial queue values of $0$. 
\end{proof} 

\begin{proof} (Theorem \ref{thm:slater}c)  Fix integers $T>0$ and $r \geq 0$.  
Similar to the proof of Lemma \ref{lem:T-slot-approx}, we have by replacing $C$ with $C(t)$ (compare
with the bound in Lemma \ref{lem:T-slot-approx}): 
\begin{eqnarray*}
\Delta_T(rT) + \sum_{\tau=rT}^{rT + T - 1} [V\hat{y}_0(\alpha(\tau), \omega(\tau)) + Vf(\bv{\gamma}(\tau)) ] \leq \\
BT + CT +  \tilde{D}T(T-1) \\
+ \sum_{\tau=rT}^{rT + T -1} [V\epsilon_V + V\hat{y}_0(\alpha^*(\tau), \omega(\tau)) + Vf(\bv{\gamma}^*(\tau)) ]  \\
+ \sum_{l=1}^LZ_l(rT) \sum_{\tau=rT}^{rT+T-1}\left[ \epsilon_Z +  \hat{y}_l(\alpha^*(\tau), \omega(\tau)) +  g_l(\bv{\gamma}^*(\tau)) \right] \\
+ \sum_{k=1}^KQ_k(rT)\sum_{\tau=rT}^{rT+T-1}\epsilon_Q  \\
+ \sum_{k=1}^KQ_k(rT)\sum_{\tau=rT}^{rT+T-1} \left[\hat{a}_k(\alpha^*(\tau), \omega(\tau)) - \hat{b}_k(\alpha^*(\tau), \omega(\tau))\right]  \\
+ \sum_{m=1}^M|H_m(rT)|\epsilon_H \\
+ \sum_{m=1}^MH_m(rT) \sum_{\tau=rT}^{rT+T-1} \left[\gamma_m^*(\tau) - \hat{x}_m(\alpha^*(\tau), \omega(\tau))\right]
\end{eqnarray*}
for any alternative decisions $\alpha^*(\tau)$, $\bv{\gamma}^*(\tau)$ over $\tau \in \{rT, \ldots, rT + T-1\}$ that
satisfy (\ref{eq:new-c4})-(\ref{eq:new-c6}).  The constant $B$ is defined according to (\ref{eq:B}) and 
$\tilde{D}$ is defined by:
\begin{equation} \label{eq:d-tilde} 
\tilde{D} \defequiv D + \frac{1}{2}\sum_{l=1}^L z_l^{diff} \epsilon_Z + \frac{1}{2}\sum_{k=1}^K q_k^{diff} \epsilon_Q + \frac{1}{2}\sum_{m=1}^Mh_m^{diff} \epsilon_H 
\end{equation} 
where $D$ is defined in (\ref{eq:D}). Note that the above bound holds \emph{deterministically} for all
possible alternative (possibly randomized) policies.  
Hence, the bound also deterministically 
holds when the right hand side is replaced by the \emph{expectation} 
over any particular randomized policy.\footnote{Formally, this uses the fact that if $b \leq \psi(\alpha_1, \ldots, \alpha_M)$ 
for all vectors $(\alpha_1, \ldots, \alpha_M) \in \script{A}$ for some function $\psi(\cdot)$, some set $\script{A}$, and some
constant $b$, 
then
$b \leq \psi(A_1, \ldots, A_M)$ for any random vector $(A_1, \ldots, A_M)$ that takes values in $\script{A}$, and hence
$b \leq \expect{\psi(A_1, \ldots, A_M)}$.}

Consider now the following randomized decisions for
$\alpha^*(t)$ and $\bv{\gamma}^*(t)$: With probability $p$ (to be defined later), for all 
slots $\tau \in \{rT, \ldots, rT + T-1\}$, choose $\alpha^*(\tau)= \alpha'_{\omega(\tau)}$ and 
$\bv{\gamma}^*(\tau) = \bv{\gamma}'(\tau)$, 
where $\alpha'_{\omega(\tau)}$ and $\bv{\gamma}'(\tau)$ are the policies from the proof of part (a) associated
with slot $\tau$.   Specifically, 
$\alpha'_{\omega(\tau)}$ satisfies Assumption A3, and $\bv{\gamma}'(\tau)$ is given by (\ref{eq:gamma-prime}). 
Else (with probability $1-p$), for all slots $\tau \in \{rT, \ldots, rT + T-1\}$ choose 
$\alpha^*(t) = \alpha''_{\omega(\tau)}$ and $\bv{\gamma}^*(\tau) = \bv{\gamma}''$, where the decisions
$\alpha''_{\omega(\tau)}$ and $\bv{\gamma}''$ 
solve (\ref{eq:frame-lookahead-problem}) and yield cost  $F_r^*$.\footnote{For simplicity, we assume
here that 
the optimum of problem (\ref{eq:frame-lookahead-problem}) is achievable by a single policy, else just take a policy
that comes within $\epsilon$ of $F_r^*$ and let $\epsilon \rightarrow 0$.} Note from our construction here 
that either all slots of the frame
use the first policy (which happens with probability $p$), or all slots of the frame use the second.
Considering the expectation of the right-hand-side under this randomized policy, we have: 
\begin{eqnarray*}
\Delta_T(rT) + \sum_{\tau=rT}^{rT + T - 1} [V\hat{y}_0(\alpha(\tau), \omega(\tau)) + Vf(\bv{\gamma}(\tau)) ] \leq \\
BT + CT +  \tilde{D}T(T-1) + VT\epsilon_V \\
+ (1-p)VTF_r^* + pT[Vy_0^{max} + Vf^{max}] \\
+ \sum_{l=1}^LZ_l(rT) \sum_{\tau=rT}^{rT+T-1}\left[ \epsilon_Z -p(\delta - (\epsilon_H+ \theta)\sum_{m=1}^M\beta_{l,m}) \right] \\
+ \sum_{k=1}^KQ_k(rT)\sum_{\tau=rT}^{rT+T-1}\left[\epsilon_Q  -p\delta\right]\\
+ \sum_{m=1}^M|H_m(rT)| \sum_{\tau=rT}^{rT+T-1} \left[\epsilon_H -  p(\epsilon_H + \theta)\right]
\end{eqnarray*}
Now choose the probability $p$ to make all of the above queueing terms non-positive, as follows: 
\[ p \defequiv \max\left[\frac{\epsilon_Z}{\delta - (\epsilon_H+\theta)\beta_{sum}}, \frac{\epsilon_Q}{\delta}, 
\frac{\epsilon_H}{\epsilon_H + \theta}\right] \]
where
\[ \beta_{sum} \defequiv \max_{l\in\{1, \ldots, L\}} \sum_{m=1}^M\beta_{l,m} \]
This is a valid probability (so that $0 \leq p \leq 1$) by definition of $\theta$ (being the solution
of (\ref{eq:theta-problem})). 
Therefore: 
\begin{eqnarray*}
\Delta_T(rT) + \sum_{\tau=rT}^{rT + T - 1} [V\hat{y}_0(\alpha(\tau), \omega(\tau)) + Vf(\bv{\gamma}(\tau)) ] \leq \\
BT + CT +  \tilde{D}T(T-1) + VT\epsilon_V \\
+ (1-p)VTF_r^* + pTV[y_0^{max} + f^{max}] 
\end{eqnarray*}
Summing the above over $r \in \{0, 1, \ldots, R-1\}$,  
using non-negativity of the Lyapunov function, all queues are initially empty,  
convexity of $f(\cdot)$, and dividing by $RTV$ yields: 
\begin{eqnarray*}
\overline{y}_0 + f(\overline{\bv{\gamma}}) &\leq& (1-p)\frac{1}{R}\sum_{r=0}^{R-1} F_r^* + \epsilon_V + p(y_0^{max} + f^{max}) \\
&& + \frac{(B + C+ \tilde{D}(T-1))}{V}
\end{eqnarray*}
Finally, we have: 
\begin{eqnarray*}
f(\overline{\bv{x}}) &\leq& f(\overline{\bv{\gamma}}) + \sum_{m=1}^M\nu_m|\overline{x}_m - \overline{\gamma}_m| \\
&\leq&  f(\overline{\bv{\gamma}}) + \sum_{m=1}^M \frac{\nu_mH_m^{max}}{RT} \\
&\leq& f(\overline{\bv{\gamma}}) + \sum_{m=1}^M \frac{\nu_m VC_3}{\theta RT}
\end{eqnarray*}
This proves Theorem \ref{thm:slater}c. 
\end{proof} 

\section*{Appendix E --- The $B$ and $D$ Constants}

Values of $B$ and $D$ that satisfy (\ref{eq:B}) and (\ref{eq:D}) are given by: 
\begin{eqnarray}
B &=&    \frac{1}{2}\sum_{l=1}^L (z_l^{diff})^2  + \frac{1}{2}\sum_{m=1}^M (h_m^{diff})^2 \nonumber \\
 && +  \frac{1}{2}\sum_{k=1}^K[(b_k^{max})^2 + (a_k^{max})^2 ] \label{eq:alwaysB} \\
D &=&  \frac{1}{2}\sum_{l=1}^L (z_l^{diff})^2  + \frac{1}{2}\sum_{m=1}^M (h_m^{diff})^2 \nonumber \\
&& + \frac{1}{2}\sum_{k=1}^K(q_k^{diff})^2   \label{eq:alwaysD}
 \end{eqnarray} 
 where constants  $z_l^{diff}$, $h_m^{diff}$, $q_k^{diff}$ are defined after equation (\ref{eq:D}).

 \emph{The Internet Model:} For the internet model of Section \ref{section:internet}, there are no 
 queues $Q_k(t)$ and so we have $q_k^{diff} = a_k^{max} = b_k^{max} = 0$.  We further have $h_m^{diff} = A_m^{max}$, and 
 \[ z_l^{diff} = \max\left[C_l^{max}, \sum_{m=1}^M 1_{l,m}A_m^{max}\right] \]
 where $1_{l,m}$ is equal to $1$ if it is possible for session $m$ to ever be routed over link $l$, and zero else. 
 Using this value of $z_l^{diff}$ in (\ref{eq:alwaysB}) and (\ref{eq:alwaysD}),  the values of $B$ and $D$ 
 for this internet 
 context, and in particular for the utility bound (\ref{eq:internet-util}), 
 are: 
 \begin{eqnarray*}
  B = D   =  \frac{1}{2}\sum_{l=1}^L (z_l^{diff})^2 + \frac{1}{2}\sum_{m=1}^M(A_m^{max})^2     
  \end{eqnarray*}

\emph{The Dynamic Queueing Network Model:} For the dynamic queueing network of Section \ref{section:network}, 
there are no $Z_l(t)$ queues and so $z_l^{diff} = 0$. Further, $h_m^{diff} = A_m^{max}$.  Because indices $k$ of queues
$Q_k(t)$ in the general framework correspond to 
indices $(n,c)$ for queues $Q_n^{(c)}(t)$ in the dynamic queueing network, 
the values of $B$ and $D$ satisfy (\ref{eq:B}) and (\ref{eq:D}) whenever the following holds
for all $t$: 
\begin{eqnarray}
B &\geq& \frac{1}{2}\sum_{m=1}^M (A_m^{max})^2 + \frac{1}{2}\sum_{n,c} [b_n^{(c)}(t)^2 + a_n^{(c)}(t)^2] \label{eq:queue-B}  \\
D&\geq&  \frac{1}{2}\sum_{n,c}\max[b_n^{(c), max}, a_n^{(c), max}]\max[b_n^{(c)}(t), a_n^{(c)}(t)]  \nonumber \\
&& + \frac{1}{2}\sum_{m=1}^M(A_m^{max})^2  \label{eq:queue-D} 
\end{eqnarray}
We use this form, rather than the more explicit form (\ref{eq:alwaysB}), (\ref{eq:alwaysD}), because this often allows
a tighter bound when we incorporate the structure of the network.  Specifically, define constants $\mu_{n}^{max, in}$, 
$\mu_n^{max, sum}$, 
$x_n^{(c), max}$ as the maximum possible sum transmission rate into node $n$, sum of transmission rates
into and out of node $n$, 
and exogenous arrivals to source node $n$ of commodity $c$, respectively, over one slot.  Note that $x_n^{(c), max}$ is given by: 
\[ x_n^{(c), max}  = \sum_{m \in \script{M}_n^{(c)}} A_m^{max} \]
Then we have  for each $n \in \{1, \ldots, N\}$:
\begin{eqnarray*}
\sum_{c} [b_n^{(c)}(t)^2 + a_n^{(c)}(t)^2]   \\
\leq \sum_{c} \left[b_n^{(c)}(t)^2 + \left(\sum_{i=1}^N \mu_{in}^{(c)}(t) + x_n^{(c), max}\right)^2\right] \\
= \sum_c \left[  b_n^{(c)}(t)^2 + \left(\sum_{i=1}^N \mu_{in}^{(c)}(t)\right)^2\right] \\
 + \sum_c \left[(x_n^{(c), max})^2 +  2\left(\sum_{i=1}^N \mu_{in}^{(c)}(t)\right)(x_n^{(c), max})\right]  \\
 \leq \left(\sum_c\left[b_n^{(c)}(t) + \left(\sum_{i=1}^N \mu_{in}^{(c)}(t)\right)\right]\right)^2 \\
 + \sum_c (x_n^{(c), max})^2  +  2\mu_{n}^{max, in}\max_{c\in\{1, \ldots, N\}} [x_n^{(c), max}] \\
 \leq (\mu_n^{max, sum})^2 \\
 + \sum_c (x_n^{(c), max})^2 + 2\mu_n^{max, in} \max_{c\in\{1, \ldots, N\}}[x_n^{(c), max}]
\end{eqnarray*}
Therefore a value of $B$ that satisfies (\ref{eq:queue-B}) is given by: 
\begin{eqnarray*}
B &=& \frac{1}{2}\sum_{n=1}^N [(\mu_n^{max, sum})^2 + \sum_c(x_n^{(c), max})^2] \\
&& + \sum_{n=1}^N \mu_n^{max, in}\max_{c\in\{1, \ldots, N\}}  [x_n^{(c), max}] \\
&& + \frac{1}{2}\sum_{m=1}^M (A_m^{max})^2
\end{eqnarray*}

Finally, define $e_n$ as follows: 
\[ e_n \defequiv \max_{c \in \{1, \ldots, N\}} \max[b_n^{(c), max}, a_n^{(c), max}] \]
Then for each $n \in \{1, \ldots, N\}$ we have: 
\begin{eqnarray*}
\sum_{c} \max[b_n^{(c), max}, a_n^{(c), max}]\max[b_n^{(c)}(t), a_n^{(c)}(t)] \\
\leq  e_n \sum_{c} \max[b_n^{(c)}(t), a_n^{(c)}(t)] \\
\leq e_n \sum_c [b_n^{(c)}(t) + a_n^{(c)}(t)] \\
\leq e_n[\mu_{n}^{max, sum} + \sum_c x_n^{(c), max}]
\end{eqnarray*}
Therefore a value of $D$ that satisfies (\ref{eq:queue-D}) is: 
\begin{eqnarray*}
D &=& \frac{1}{2}\sum_{m=1}^M(A_m^{max})^2  \\
&& + \frac{1}{2}\sum_{n=1}^N e_n[\mu_n^{max, sum} + \sum_c x_n^{(c), max}]
\end{eqnarray*}

For example, consider a wireless network where data is measured in integer units of packets (assumed to have a fixed
length).  Suppose that at most one packet can be transmitted or received per node per slot, and that a packet cannot
be transmitted and received on the same slot at the same node.  Then 
we have $\mu_n^{max, sum} = \mu_n^{max, in} = 1$.  Further, 
suppose 
there is at most one source at any given node (so that $M \leq N$), and no source can admit more than 1 packet per slot.
Then $\sum_c x_n^{(c), max} = 1$ if node $n$ is a source, and zero else, and $e_n = 2$ if node $n$ is a source, 
and $1$ else.   There are $M$ source nodes and $N-M$ non-source nodes, and so
$B$ and $D$ are: 
\begin{eqnarray*}
B = D = (N + 4M)/2 
\end{eqnarray*}

\bibliographystyle{unsrt}
\bibliography{../../latex-mit/bibliography/refs}

\begin{thebibliography}{10}

\bibitem{now}
L.~Georgiadis, M.~J. Neely, and L.~Tassiulas.
\newblock Resource allocation and cross-layer control in wireless networks.
\newblock {\em Foundations and Trends in Networking}, vol. 1, no. 1, pp. 1-149,
  2006.

\bibitem{neely-fairness-infocom05}
M.~J. Neely, E.~Modiano, and C.~Li.
\newblock Fairness and optimal stochastic control for heterogeneous networks.
\newblock {\em Proc. IEEE INFOCOM}, March 2005.

\bibitem{neely-fairness-ton}
M.~J. Neely, E.~Modiano, and C.~Li.
\newblock Fairness and optimal stochastic control for heterogeneous networks.
\newblock {\em IEEE/ACM Transactions on Networking}, vol. 16, no. 2, pp.
  396-409, April 2008.

\bibitem{neely-thesis}
M.~J. Neely.
\newblock {\em Dynamic Power Allocation and Routing for Satellite and Wireless
  Networks with Time Varying Channels}.
\newblock PhD thesis, Massachusetts Institute of Technology, LIDS, 2003.

\bibitem{tass-one-hop}
L.~Tassiulas.
\newblock Scheduling and performance limits of networks with constantly
  changing topology.
\newblock {\em IEEE Trans. on Inf. Theory}, May 1997.

\bibitem{neely-power-network-jsac}
M.~J. Neely, E.~Modiano, and C.~E Rohrs.
\newblock Dynamic power allocation and routing for time varying wireless
  networks.
\newblock {\em IEEE Journal on Selected Areas in Communications}, vol. 23, no.
  1, pp. 89-103, January 2005.

\bibitem{neely-maximal-bursty-ton}
M.~J. Neely.
\newblock Delay analysis for maximal scheduling with flow control in wireless
  networks with bursty traffic.
\newblock {\em IEEE Transactions on Networking}, vol. 17, no. 4, pp. 1146-1159,
  August 2009.

\bibitem{neely-switch}
M.~J. Neely, E.~Modiano, and C.~E. Rohrs.
\newblock Tradeoffs in delay guarantees and computation complexity for $n
  \times n$ packet switches.
\newblock {\em Proc. of Conf. on Information Sciences and Systems (CISS),
  Princeton}, March 2002.

\bibitem{neely-mesh}
M.~J. Neely and R.~Urgaonkar.
\newblock Cross layer adaptive control for wireless mesh networks.
\newblock {\em Ad Hoc Networks (Elsevier)}, vol. 5, no. 6, pp. 719-743, August
  2007.

\bibitem{nequit-paper}
J.~Andrews, S.~Shakkottai, R.~Heath, N.~Jindal, M.~Haenggi, R.~Berry, D.~Guo,
  M.~Neely, S.~Weber, S.~Jafar, and A.~Yener.
\newblock Rethinking information theory for mobile ad hoc networks.
\newblock {\em IEEE Communications Magazine}, vol. 46, no. 12, pp. 94-101, Dec.
  2008.

\bibitem{lempel-ziv}
J.~Ziv and A.~Lempel.
\newblock A universal algorithm for sequential data compression.
\newblock {\em IEEE Transactions on Information Theory}, IT-23, no. 3, pp.
  337-343, May 1977.

\bibitem{cover-universal}
T.~M. Cover.
\newblock Universal portfolios.
\newblock {\em Mathematical Finance}, vol. 1, no. 1, pp. 1-29, Jan. 1991.

\bibitem{universal-stock2}
T.~M. Cover and E.~Ordentlich.
\newblock Universal portfolios with side information.
\newblock {\em IEEE Transactions on Information Theory}, vol. 42, no. 2, 1996.

\bibitem{merhav-universal}
N.~Merhav and M.~Feder.
\newblock Universal schemes for sequential decision from individual data
  sequences.
\newblock {\em IEEE Transactions on Information Theory}, vol. 39, no. 4, pp.
  1280-1292, July 1993.

\bibitem{neely-stock-arxiv}
M.~J. Neely.
\newblock Stock market trading via stochastic network optimization.
\newblock {\em ArXiv Technical Report}, arXiv:0909.3891v1, Sept. 2009.

\bibitem{competitive-ratio-garay}
J.~A. Garay and I.~S. Gopal.
\newblock Call preemption in communication networks.
\newblock {\em Proc. IEEE INFOCOM}, vol. 44, pp. 1043-1050, Florence, Italy,
  1992.

\bibitem{plotkin-competitive}
S.~Plotkin.
\newblock Competitive routing of virtual circuits in atm networks.
\newblock {\em IEEE Journal on Selected Areas in Communications}, vol. 13, no.
  6, pp. 1128-1136, Aug. 1995.

\bibitem{lin-shroff-large-N-energy}
L.~Lin, N.~B. Shroff, and R.~Srikant.
\newblock Asymptotically optimal power-aware routing for multihop wireless
  networks with renewable energy sources.
\newblock {\em Proc. IEEE INFOCOM}, March 2005.

\bibitem{srikant-universal}
J.~J. Jaramillo and R.~Srikant.
\newblock Admission control and routing in multi-hop wireless networks.
\newblock {\em Proc. IEEE Conf. on Decision and Control (CDC)}, Dec. 2008.

\bibitem{andrews-max-profit}
M.~Andrews.
\newblock Maximizing profit in overloaded networks.
\newblock {\em Proc. IEEE INFOCOM}, March 2005.

\bibitem{stolyar-greedy}
A.~Stolyar.
\newblock Maximizing queueing network utility subject to stability: Greedy
  primal-dual algorithm.
\newblock {\em Queueing Systems}, vol. 50, pp. 401-457, 2005.

\bibitem{stolyar-gpd-gen}
A.~Stolyar.
\newblock Greedy primal-dual algorithm for dynamic resource allocation in
  complex networks.
\newblock {\em Queueing Systems}, vol. 54, pp. 203-220, 2006.

\bibitem{kelly-shadowprice}
F.P. Kelly, A.Maulloo, and D.~Tan.
\newblock Rate control for communication networks: Shadow prices, proportional
  fairness, and stability.
\newblock {\em Journ. of the Operational Res. Society}, 49, p.237-252, 1998.

\bibitem{low-flow-control}
S.~H. Low and D.~E. Lapsley.
\newblock Optimization flow control, i: Basic algorithm and convergence.
\newblock {\em IEEE/ACM Transactions on Networking}, vol. 7(6): 861-75, Dec.
  1999.

\bibitem{chiang-layering-decomposition}
M.~Chiang, S.~H. Low, A.~R. Calderbank, and J.~C. Doyle.
\newblock Layering as optimization decomposition: A mathematical theory of
  network architectures.
\newblock {\em Proceedings of the IEEE}, vol. 95, no. 1, Jan. 2007.

\bibitem{leonardi-SP-routing}
E.~Leonardi, M.~Mellia, M.~A. Marsan, and F.~Neri.
\newblock Optimal scheduling and routing for maximizing network throughput.
\newblock {\em IEEE/ACM Transactions on Networking}, vol. 15, no. 6, Dec. 2007.

\bibitem{chiang-delays-flownets}
Y.~Li, A.~Papachristodoulou, and M.~Chiang.
\newblock Stability of congestion control schemes with delay sensitive traffic.
\newblock {\em Proc. IEEE ACC, Seattle, WA}, June 2008.

\bibitem{pricing-congested-links}
J.~K. MacKie-Mason and H.~R. Varian.
\newblock Pricing congestible network resources.
\newblock {\em IEEE Journal on Selected Areas in Communications}, vol. 13, no.
  7, September 1995.

\bibitem{bertsekas-data-nets}
D.~P. Bertsekas and R.~Gallager.
\newblock {\em Data Networks}.
\newblock New Jersey: Prentice-Hall, Inc., 1992.

\bibitem{neely-convex-it}
M.~J. Neely and E.~Modiano.
\newblock Convexity in queues with general inputs.
\newblock {\em IEEE Transactions on Information Theory}, vol. 51, no. 2, pp.
  706-714, Feb. 2005.

\bibitem{lin-shroff-cdc04}
X.~Lin and N.~B. Shroff.
\newblock Joint rate control and scheduling in multihop wireless networks.
\newblock {\em Proc. of 43rd IEEE Conf. on Decision and Control, Paradise
  Island, Bahamas}, Dec. 2004.

\bibitem{neely-energy-it}
M.~J. Neely.
\newblock Energy optimal control for time varying wireless networks.
\newblock {\em IEEE Transactions on Information Theory}, vol. 52, no. 7, pp.
  2915-2934, July 2006.

\bibitem{bertsekas-nonlinear}
D.~P. Bertsekas.
\newblock {\em Nonlinear Programming}.
\newblock Athena Scientific, Belmont, MA, 1995.

\bibitem{neely-delay-based-infocom2010}
M.~J. Neely.
\newblock Delay-based network utility maximization.
\newblock {\em Proc. IEEE INFOCOM}, 2010.

\end{thebibliography}
\end{document}